\documentclass[leqno,11pt]{article}

\usepackage{amssymb}

\usepackage{euscript}
\evensidemargin\oddsidemargin

\renewcommand{\theequation}{\thesection.\arabic{equation}}  
\newtheorem{theorem}{Theorem}[section] 
\newtheorem{lemma}{Lemma}[section]
\newtheorem{fact}{Fact}
\newtheorem{Proposition}{Proposition}[section]
\newcommand{\eqnsection}{
\renewcommand{\theequation}{\thesection.\arabic{equation}}
    \makeatletter
    \csname  @addtoreset\endcsname{equation}{section}
    \makeatother}
\eqnsection

\def\r{{\mathbb R}}
\def\e{{\mathbb E}}
\def\b{{\mathbb B}}
\def\p{{\mathbb P}}
\def\z{{\mathbb Z}}
\def\ee{\mathrm{e}} 
\def\d{\, \mathrm{d}}
\def\law{{\buildrel \mathrm{law} \over =}}
\def\ddef{{\buildrel \mathrm{def} \over =}}
\def\qed{\hfill$\square$}
\makeatletter
\newcommand{\silentfootnote}[1]
   {\gdef\@thefnmark{~}\@footnotetext{#1}}
\makeatother
\title{\bf Annealed tail estimates for a Brownian motion in a drifted Brownian
potential} \author{Marina Talet}

\date{}

\begin{document}
\maketitle
\silentfootnote{Key words and phrases. Large deviation, Brownian
motion in a random potential, Lamperti's representation, drifted
Brownian motion, Bessel process.}  \silentfootnote{AMS (1991) subject classifications:
60K37, 60J55.}

\begin{abstract}
We study Brownian motion in a drifted
Brownian potential. Kawazu and Tanaka \cite{Kawazu-Tanaka2} exhibited two
speed regimes for this process, depending on the drift. They supplemented these laws of large numbers by central limit 
theorems, which were recently completed
by Hu, Shi and Yor \cite{HSY} using stochastic calculus. We studied large
deviations \cite{Marina}, showing among other results that the rate function in the annealed setting, that is after averaging over the potential, has a
flat piece in
the ballistic regime. In this paper, we focus on this subexponential regime,
proving that the probability of deviating below the almost sure speed has a
polynomial rate of decay, and computing the exponent in this power law. This
provides the continuous-time analogue of what Dembo, Peres and Zeitouni proved
for the transient random walk in random environment \cite{DPZ}. Our method
takes a completely different route, making use of Lamperti's representation
together with an iteration scheme.

\end{abstract}

%\vfill\eject

\section{Introduction}
   \label{s:intro}

Let $\omega=\{\omega_i\}_{i\in \z}$ be a collection of i.i.d. $(0,1)$-valued
random variables, serving as an environment, and define a conditional Markov 
chain on the integers, $\{S_n\}_{n\ge 0}$, by $S_0=0$ and
%taking values in $(0,1)$, such that the support of the
%law of $\omega_0$ is contained in $(0,1)$. 

$$
\p\left( S_{n+1}=y \, | \, S_n=x, \; \{ \omega_i\}_{i\in \z}
\right)
= \cases{ \omega_x &if $y=x+1$, \cr
         1-\omega_x &if $y=x-1$, \cr
         0 &otherwise. \cr}
$$

\noindent The process $\{S_n\}_{n\ge 0}$ is called a random walk
in a random environment (hereafter abbreviated RWRE).

\medskip

Solomon \cite{Solomon} completely solved the transience/recurrence
problem for $\{S_n\}_{n\ge 0}$, and determined furthermore the
speed of the walk. In particular, setting $\rho\; \ddef\;
(1-\omega_0)/\omega_0$, he proved that if $\e(\rho)<1$,
then, almost surely, $\lim_{n\to \infty} S_n /n= (1-\e \rho)/(1+\e
\rho)\stackrel{def}{=}$v (with the case $\e(1/\rho)<1$ then following by
reflection), and that otherwise $\lim_{n\to \infty}
S_n /n = 0$ almost surely. These laws of large numbers
were later developed into central limit theorems by Kesten, Kozlov and Spitzer
\cite{Kesten-Kozlov-Spitzer}.

\medskip

Large deviations for $\{S_n/n\}_{n> 0}$ were investigated
by several authors, both under the conditional probability given $\omega$, the
so-called {\it quenched} probability $P^{\omega}$, and the {\it annealed} one
$\p$, that given after averaging over the environment $\omega$. We refer to Greven and den Hollander \cite{gdh}, Gantert and Zeitouni \cite{Gantert-Zeitouni}
and Comets, Gantert and Zeitouni \cite{Comets-Gantert-Zeitouni} for insightful
overviews.  

\bigskip

Dembo, Peres and Zeitouni \cite{DPZ} studied the probability with which the
walk deviates from its limiting speed when this speed is nonzero. Assuming
$\e(\rho)<1$ and $\p(\omega_0<1/2)>0$, let us write $s$ for the unique $s>1$ such that \begin{equation}
    \e(\rho^s)=1.
    \label{s}
\end{equation}

\smallskip

\noindent They proved that 

\bigskip

\noindent {\bf Theorem A} (\cite{DPZ}) {\it For any
open $G\subset$ (0,}v{\it) which is separated from} v,  
$$
\lim_{n\to \infty} \frac{1}{\log n}\log \p\left(\frac{S_n}{n} \in G\right)=1-s.$$

\bigskip

\noindent This result was the starting point for our work.

\bigskip

\noindent In the present paper, we are interested in the continuous-time analogue of
RWREs, the so-called Brownian motion in random potential $W=\{W(x)\}_{x\in
  \r}$. This should be a solution of the {\it formal} stochastic differential equation
$$
    \cases{ \d X(t) = \d \alpha(t) - {1\over 2}\, W'(X(t)) \d t, \cr 
    X(0)=0, \cr } 
 $$

\medskip

\noindent where the potential $W$ is defined by
$$
    W(x)= B(x) - {\kappa\over 2} x, \qquad x\in \r, \qquad \kappa \in \r,
$$   

\noindent with $\{B(x)\}_{x\in \r}$ is a one-dimensional two-sided Brownian
motion defined on $\r$ starting from zero and $\{\alpha(x)\}_{x\ge 0}$ a
standard Brownian motion such that $\alpha(0)=0$, independent of $W$ (or
equivalently of $B$).

\medskip

One way of defining a ``formal solution'' is this: since the Brownian motion is almost surely
nowhere differentiable, one defines the process $X$ through its conditional generator given $W$, 
$$ 
{\cal L}_{W}\; \ddef\; {1\over 2} \, \ee^{W(x)} {\! \d\over \! \d x} \left( \ee^{-W(x)}
{\! \d \over \! \d x} \right) .
$$

\noindent Since we are dealing with one-dimensional diffusions, there is a second
approach to defining $X$, which we shall adopt. The martingale representation
for diffusions tells us that, at fixed environment, that is to say for each realization of
the environment $W$, the image of $X$ under its scale function, which is a
continuous martingale, can be represented as a time-changed Brownian
motion. Namely, 
\begin{equation}
    X(t) = S^{-1} ({\cal B}(T^{-1}(t))), \qquad t\ge 0, 
    \label{Brox}
\end{equation} 

\noindent where $\{{\cal B}(t)\}_{t\ge 0}$ is a standard Brownian motion starting from 0, independent of $W$, with the scale function $S$ and random clock
$T$ defined by
\begin{eqnarray}
    S(x) 
 &\; \ddef \;& \int_0^x \ee^{W(y)} \d y, \qquad x\in \r, 
    \label{S} \\
    T(t)
 &\; \ddef\; & \int_0^t \exp \left( -2W(S^{-1}({\cal B}(u)))\right) \d u, \qquad
    t\ge 0,
    \label{T} 
\end{eqnarray}

\noindent where $S^{-1}$ and $T^{-1}$ denote the respective inverse functions 
of $S$ and $T$.\\

In the {\it quenched} situation, i.e. at fixed environment, $X$ is Markov. We
denote its law by $P^W$ and the Wiener measure by $Q$. Averaging
$P^W$ over $Q$ gives birth to a new probability $\p$, called the {\it
  annealed} probability. Note
that, under $\p$, the process $X$ is not necessarily Markov.

\bigskip

Brox \cite{Brox} was the first to study such processes. He proved that for $\kappa=0$, in which case the diffusion is recurrent, 
the motion is extremely slow, as then $X(t)$ is of order $\log^2 t$ for large $t$, in this way differing markedly from the diffusive behavior of Brownian motion.

\medskip

For $\kappa \not=0$, $X$ is transient to the left or the right depending on the
sign of $\kappa$, which is just ``space reversal invariance''. 
Kawazu and Tanaka \cite{Kawazu-Tanaka2} computed its almost-sure speed; assuming $\kappa>0$, 

$$
    \lim_{t\to \infty} {X(t)\over t}= v_\kappa\; \ddef\; {(\kappa-1)^+\over
    4}, \qquad \p- \hbox{\rm  a.s.}
$$    

\noindent  These are
continuous-time analogues of Solomon's aforementioned laws of large
numbers for RWRE. The corresponding
central limit theorems were established by Kawazu and Tanaka
\cite{Kawazu-Tanaka3} using Krein's spectral theory, and both recovered and
completed by Hu, Shi and Yor \cite{HSY} using stochastic calculus. We proved in \cite{Marina} that the family
of distributions of $\{X(t)/t\}_{t> 0}$ satisfies a Large Deviation
Principle in both the quenched and the annealed frameworks. We note that
$\kappa$ plays the role of $s$ (defined for the RWRE in (\ref{s})) for these
laws of large numbers, for central limit theorems as well as for the results of the present paper.

\bigskip

As in the discrete case, we are interested in
the probability with which $X$ deviates from its limiting speed
$v_\kappa$ when $v_\kappa\not= 0$, in the annealed setting. By symmetry, we
only have to deal with $\kappa>1$ in which case, according to \cite{Marina},
the annealed rate function has a ``flat'' piece, by which one means that
$$
\lim_{t\to \infty} \frac{1}{t} \log \p\left(\frac{X(t)}{t}\in G\right)=0,
$$

\medskip

\noindent for any open set $G \subset (0, v_\kappa)$ which is
separated from $v_\kappa$. This tells us that the probability for $X(t)/t$ to
deviate below the typical velocity decays subexponentially fast to zero. 

\smallskip

A natural question arises : how fast exactly? And can a result similar to Theorem A be
obtained for $X$? The answer is provided by:
\bigskip
\begin{theorem}
 \label{t:main}
   Let $\kappa>1$. For any open $G\subset (0, v_\kappa)$ which
   is separated from $v_\kappa$, 
   $$ \lim_{t\to \infty} \frac{1}{\log t} \log \p\left(\frac{X(t)}{t} \in G\right) =1-\kappa.$$
      \end{theorem}
 
\bigskip

\noindent Again, as in \cite{Marina}, following the approach of \cite{Comets-Gantert-Zeitouni}, we set 
\begin{equation}
    H(r) \; \ddef\; \inf \{ t>0: \; X(t) > r\}, \qquad r>0.
    \label{H}
\end{equation}

\noindent Making use of the ``natural duality'' between the diffusion $X$ and its first hitting time process $H$, that is to say
$$
\p\left( \frac{X(t)}{t} \sim \; v \right)\: \approx\:  \p\left(\frac{H(tv)}{t}
  \: \sim  \frac{1}{v}\right),$$

\noindent the proof of Theorem \ref{t:main} reduces to showing 

\medskip

\begin{theorem}
 \label{t:H}
    Let $\kappa>1$. For any $u>v_\kappa^{-1}=4/(\kappa-1)$, 
    \begin{equation}
       \limsup_{r\to \infty} \frac{1}{\log r} \log \p(H(r)>ur)  \le 
       1-\kappa,
       \label{t:H-ub}
    \end{equation}
    and for any $u>0$, 
    \begin{equation}
       \liminf_{r\to \infty} \frac{1}{\log r} \log \p(H(r)>ur) \ge 
       1-\kappa. 
       \label{t:H-lb}
    \end{equation}
\end{theorem}

\bigskip

Differently from \cite{DPZ}, our proof hinges upon stochastic calculus
  techniques. A key role is played by Lamperti's representation which relates the potential $W$ to a Bessel
  process. This fact enabled Hu, Shi and Yor to derive central limit theorems
  for the model we are studying, in \cite{HSY}. Here we are interested in deviation estimates, hence in the rate at which
various random variables involved in \cite{HSY} converge. This leads to the delicate probability
estimates of Sections 5 and 6.

\smallskip

We note from a glance at both \cite{Kawazu-Tanaka3} and \cite{HSY} that
in the case where $1<\kappa<2$, 
$$ \frac{1}{r^{1/\kappa}} \left(H(r)-\frac{4}{\kappa-1}\: r\right)
  \stackrel{law}{\longrightarrow} {\mathrm {\; \: stable \; variable}},$$

\noindent where $\stackrel{law}{\longrightarrow}$ denotes convergence in
  distribution. In this case, as we proved in \cite{Marina3}, the main contribution
  to the polynomial rate of decay of $\p(H(r)>ur)$, ($\sim r^{1-\kappa}$),
  stems from the limiting stable law in this regime. 

\medskip

 To prove Theorem \ref{t:H}, we introduce an iteration technique which (as far as we know)
  is new and may prove to be of use elsewhere. 

%We apply successively
 % integration by parts followed by a time change. (See (\ref{mimi3}) and (\ref{bonn})). 

\medskip
 Solely using integration by parts followed by a time change (see (\ref{mimi3}) and (\ref{bonn})), together with results on Bessel 
and stable processes, our iteration 
blends very naturally with the techniques used in this paper. It also offers an alternative strategy that 
circumvents dealing with special functions while solving a Sturm-Liouville equation; see Section \ref{s:Appendix}.

%circumventing heavy computations 
%delaing with special functions while solving a Sturm-Liouville equation, as is shown in Section \ref{s:Appendix}. 

\bigskip

\noindent The outline of the paper is as follows. In Section
\ref{s:preliminaries}, we list a collection of known results on Brownian local
times, Bessel and Jacobi processes. In Section \ref{s:Two lemmatas}, we state
and prove three lemmas which will be of frequent use in the proof of our tail
estimates. The main result of this paper is Theorem \ref{t:main}. As in the 
discrete case, proving Theorem \ref{t:main} reduces to proving Theorem \ref{t:H} for the first
hitting time process. This step is justified in Section \ref{s:main}, and
the rest of the paper is devoted to the proof of Theorem \ref{t:H}. This
proceeds in one further step: Theorem
\ref{p:derniere}, stated in Section \ref{s:Second reduction}. We prove Theorem \ref{p:derniere} in Section 
\ref{s:six} and Section \ref{half} by means of a key estimate stated and proved in Section 
\ref{s:key}. Theorem \ref{p:derniere} implies Theorem \ref{t:H}; this is proved in Section \ref{s:Second reduction}. 
And  the last section is devoted to solving a Sturm-Liouville equation, providing an alternative method to our iteration technique. 

\medskip

\noindent{\bf Notation:} Throughout the sequel, $Q$ will denote the Wiener
measure, $E_Q$ the expectation with respect to $Q$, $P^W_x$ and 
$\p_x\stackrel{def}{=}E_Q(P^W_x(.))$ the quenched and
annealed laws when $X(0)=x$, and $E^W_{x}$ and $\e_x$ the
expectations w.r.t. $P^W_x$ and $\p_x$ respectively. For notational convenience, $P^W$, $\p$, $E^W$ and $\e$ stand for $P^{W}_{0}$, $\p_0$,
$E^W_0$ and $\e_0$. We sometimes drop the subscript $x$ in $\p_x$ or in $\e_x$, for $x\neq0$, when no confusion can arise.  

\bigskip

\noindent Unless stated otherwise, {\bf it is assumed that
$\kappa>1$}.

\section{Preliminaries}
\label{s:preliminaries}

In this section, we summarize a collection of known results
which will be useful in the rest of the paper. These results
concern Brownian local times, Bessel and Jacobi processes as well as 
Lamperti's representation for exponentials of drifted Brownian motions.\\

Let $\{\gamma(t)\}_{t\ge 0}$ be a standard Brownian motion. A well-known
theorem of Trotter \cite{Trotter} confirms the existence of a jointly continuous version of the local time process
$\{L_t^x (\gamma)\}_{t\ge 0, \, x\in \r}$ as the density of occupation
time: for any bounded Borel function $f$, 
\begin{equation}
    \int_0^t f(\gamma(s))\d s= \int_\r f(x) L_t^x (\gamma) \d x . 
    \label{occupation}
\end{equation}

\noindent Let 
\begin{eqnarray}
\sigma_{\gamma}(r) &\; \ddef \; & \inf\{ t>0: \; \gamma(t)>r\}, \qquad r>0, \label{sigma} \\   \tau_{\gamma}(r) &\; \ddef \; &  \inf\{ t>0: \; L_t^0 (\gamma) >r\}, 
\qquad r>0,
\label{tau}
\end{eqnarray}
denote the first hitting time of $\gamma$ and its inverse local time at 0 respectively. We shall drop the subscript $\gamma$ when no confusion arises.
\smallskip

As is shown by the Ray--Knight theorems (see Revuz and Yor \cite{R-Y},
Chap XI),
Brownian local times at these hitting times are nice diffusion processes,
known as Bessel processes.

\bigskip

\noindent {\bf Definition.} {\it A squared Bessel process $\{R^2(t), t\ge 0, \p\}$ of
  dimension $d$ and started at $r^2$, is the solution of the stochastic differential equation
\begin{equation}
\label{bessel:def} \d R^2 (t) = 2 R(t) \d\gamma(t) + d \: \d t,
\end{equation}
with $R^2 (0)=r^2$ and $\gamma$ a $\p$-Brownian motion. A Bessel process of dimension $d$, started at $r$, is $\{R(t), t\ge 0
, \p\}$ with $R(t)\geq 0$ and $R(0)=r$. 
}

\bigskip

\noindent We recall seven facts from the literature.

\medskip

\begin{fact} (Ray-Knight theorems)  
\smallskip

\noindent {\bf First:} The process $\{L_{\sigma(1)}^{1-t}\}_{t\ge 0}$ is a
squared Bessel process started at 0, of dimension 2 for $0\le t \le 1$ and of
dimension 0 for $t\ge 1$. (Here $L_{\sigma(1)}^0$ is an exponential random variable of mean 2.)
 
\medskip

\noindent {\bf Second:} The process $\{L_{\tau(1)}^t\}_{t\ge 0}$ is a squared
Bessel process of dimension 0, started at 1.
\end{fact}

%\medskip

\begin{fact} (Revuz and Yor, \cite{R-Y}, p 430) Let $\{R (t)\}_{t
\ge 0}$ be a Bessel process of dimension $d>2$, starting from $x>0$. Then
\begin{eqnarray}
\label{lam1}\lim_{r \rightarrow \infty} \frac{1}{\log r} \int_0^r \frac{\d u}{R^2 (u)}=\frac{1}{d-2},\qquad  \p \: a.s. 
\end{eqnarray}
%Furthermore, if $d\ge 3$, we have, for all $p>0$, 
%\begin{eqnarray}
%\label{lam2} \sup_{r>1} \frac{1}{\log^p r} \e\left(\left(\int_0^r \frac{\d u}{R^2 (u)}\right)^p\right)&<& \infty.
%\end{eqnarray}
\end{fact}

\bigskip

\noindent The following result was first proved by Dufresne \cite{Dufresne}
using direct computations. We learned it from Yor \cite{Yor}. 

\bigskip

\begin{fact} (Dufresne) {\it Let $\kappa>0$. The law of the almost sure
    random variable $S(\infty)$ is, up to a constant, the inverse
    of a Gamma distribution. More precisely,}
\begin{equation}
    \p\left( S(\infty)\in \! \d x\right) = {2^\kappa \over
    \Gamma(\kappa)} \ee^{-2/x} x^{-(\kappa+1)} \d x, \; \; {\mathrm for}\; \qquad x>0.
    \label{Gamma}
\end{equation}
\end{fact}

\bigskip

\noindent A powerful tool in the study of
exponential functionals of drifted Brownian motions, and more generally of L\'evy
processes, is Lamperti's representation.

\medskip

\begin{fact} (Lamperti, \cite{Lamperti}) Let $\zeta \in \r$. There exists $\{R(t)\}_{t\ge 0}$, a Bessel process of dimension
$(2+2\zeta)$ starting from $2$, such that
\begin{equation}
    \ee^{B(t)+\zeta t/2}= \frac{1}{4} R^2 \left(\int_0^t e^{B(y)+\zeta y/2} \d y \right), \qquad
    t\ge 0.
    \label{Lamperti}
\end{equation} 
\end{fact}

\smallskip

\noindent In particular, taking $\zeta=-\kappa$, $S(\infty)$ appears as the first hitting time of 0 by $R$. Recall that in this case $R$ is a Bessel process of dimension $2-2\kappa$. 
\bigskip

\begin{fact} (Getoor and Sharpe, \cite{Getoor-Sharpe}) For all $z \ge 0$, and for all $u\geq0$ such that $2uz<1$, we have
\begin{equation}
\label{bass-griffin}
\e\left(\exp(u L_{\tau (1)}^z) \right)= \exp\left(\frac{u}{1-2uz}\right).
\end{equation}
\end{fact}

\medskip

\noindent Moreover,

\begin{fact} (Biane and Yor, \cite{Biane-Yor}) For any $\lambda>0$, $0<p<1$, 
\begin{eqnarray}\label{bessstab1}\int_0^{\infty} x^{1/p-2} 
L_{\tau (\lambda)}^x \: dx&\; \law \;& 
2 p^{2-2/p} \psi(p) \lambda^{1/p}
 {\cal S}_{p},\\
\int_0^1 \frac{L_{\tau (1)}^x -1}{x} \: dx + \int_1^{\infty} \frac{L_{\tau(1)}^x}{x} dx &\; \law \;& 2{\bar {\gamma}} +\log \frac{\pi}{4} + \frac{\pi}{2} {\cal C}_1,\label{cauchy}
\end{eqnarray}

\medskip

\noindent where $$\psi(p)=\left(\frac{\pi p}{4 \Gamma^2(p) 
\sin(\pi p/2)}\right)^{1/p},$$

\noindent with ${\bar {\gamma}}$ denoting Euler's constant, $\Gamma$ the usual
gamma function, ${\cal S}_{p}$ a completely
asymmetric stable variable of index $p$ and ${\cal C}_1$ a completely
asymmetric Cauchy variable of index 1. The laws of ${\cal S}_p$ and ${\cal C}_1$ are characterized by
\begin{eqnarray*}
\e\left(\ee^{it {\cal
      S}_{p}}\right)&=&\exp\left(-|t|^{p}\left(1-i \: {\rm sgn} (t) \tan \frac{\pi p}{2}\right)\right),\\ 
\e\left(\ee^{it {\cal C}_{1}}\right)&=& \exp\left(-|t|-it\frac{2}{\pi} \log|t|\right).
\end{eqnarray*}\end{fact}
\noindent For a reference on stable laws, see e.g. \cite{Bingham}, p 347.

\bigskip

\noindent Before stating the final result, let us recall the definition of a Jacobi
process. See for instance \cite{K-T}.

\bigskip 

\noindent {\bf Definition.} {\it A 
Jacobi process $\{Y(t),\: t\geq 0, \p\}$ 
of dimensions $(d_1, d_2)$ starting from $y \in (0,1)$ is the solution of the stochastic differential equation
\begin{equation} \label{jacobieds}\d Y(t)=2\sqrt{Y(t)(1-Y(t))}\; \d \b(t)+
(d_1-(d_1+d_2)Y(t)) \d t,\end{equation}
with $0\leq Y(t) \leq 1$, $Y(0)=y$ and with $\b$ a $\p$-Brownian motion.}

\bigskip

\noindent The following result is due to Warren and Yor \cite{Warren-Yor}; it 
relates Bessel and Jacobi processes.
\medskip

\begin{fact} (Warren and Yor) Let $\{R_1 (t)\}_{t \ge 0}$ and 
$\{R_2 (t)\}_{t \ge 0}$ be two independent Bessel processes of dimensions $d_1$ and $d_2$ respectively, with $d_1 + d_2 \ge 2$, $R_1 (0)=r_1 \ge0$ and $R_2 (0)=r_2 >0$. There exists a Jacobi process $\{Y(t)\}_{t \ge 0}$ of dimensions $(d_1, d_2)$, starting from $r_1^2/ (r_1^2 + r_2^2)$, independent of 
$\{R_1^2 (t)+R_1^2 (t)\}_{t \ge 0}$, and such that for all $t \ge 0$,  
\begin{equation}\label{warren}\frac{R_1^2(t)}{R_1^2(t)+R_2^2(t)}=Y\left(\int_0^t
    \frac{\d s}{R_1^2(s)+
R_2^2(s)}\right).\end{equation} \end{fact}

\bigskip

\noindent Let us now state and prove three lemmas which will be of constant use in what follows.

\section{Three lemmas}
\label{s:Two lemmatas}
\begin{lemma}\label{uno} 
Let $\{Z(t)\}_{t\ge 0}$ denote a Bessel process of dimension 0 started at 1. For all $v, \delta>0$ and $u\geq1$, we have
\begin{eqnarray} 
\label{bess1}\p\left(\sup_{t\geq 0} Z(t)>u\right)&=&\frac{1}{u},\\
\label{bess2}\p\left(\sup_{0\leq t\leq v} |Z(t)-1| > \delta\right) &\leq& 
4 \frac{\sqrt{(1+\delta) v}}{\delta} \exp\left(-\frac{\delta^2}{8(1+\delta)
    v}\right).\end{eqnarray}\end{lemma}

\bigskip

\noindent Set
\[\Sigma(r)\stackrel{def}{=} \int_0^r e^{-W(y)} \d y.\]
%\[\Lambda(r)\stackrel{def}{=} \int_0^r \frac{1}{R^2 (s) }\d s,\]
Intuitively speaking, $\log \Sigma(r)$ is of order $\kappa r/2$. The following
lemma gives a rigorous form to this intuition.

\begin{lemma} \label{due}
For any $\delta>0$, there exist two constants $c_1$ and $c_2$ depending on
both $\delta$ and $\kappa$, such that, for $r$ big enough, 
\begin{equation}\label{Sigma}\p\left(\left|\log \Sigma (r)-\frac{\kappa}{2} r\right|>\delta r\right)\le
c_1 e^{-c_2 r}.\end{equation}

\bigskip

Furthermore, let $\{ R(t)\}_{t\ge 0}$ denote a Bessel process of dimension $d>2$, starting
at 2. For any $\delta>0$, there
  exist two constants $c_3$ and $c_4$ depending on both $\delta$ and $d$ such that, for all $r$ big enough,
\begin{equation}\label{lam3}\p \left( \left|\frac{1}{\log r}\int_0^r \frac{\d
  s}{R^2 (s)} -\frac{1}{d-2} 
\right|> \delta \right) \le \frac{c_3}{r^{c_4}}.\end{equation}
\end{lemma}

\bigskip

\noindent The last result complements (\ref{lam1}). Still dealing with Bessel processes, the following lemma will be used 
in Section \ref{s:End}.  

\medskip
\begin{lemma}
 \label{l:Bessel-ub}  
Let $\{ R(t)\}_{t\ge 0}$ denote a Bessel process of dimension $d>2$, starting
at $x>0$. Let $\, a\ge 0$ and $b>2a+2$. For any
    $p<(d-2)/(b-2)$, we have
    $$
      \e \left( \, \int_0^\infty {s^a \over R^b(s)} \d s
      \right)^{\! p} < \infty.    $$
\end{lemma}

%\noindent with $\p_x$ and $\e_x$ denoting probability and expectation w.r.t $\p_x$ when the process starts at $x$.  

\bigskip

\noindent Let us prove the aforestated lemmas; we begin with the

\bigskip

\noindent{\bf Proof of Lemma \ref{uno}.} Recall from (\ref{bessel:def}) with
$d=0$ that the process $Z$ solves
\begin{equation}
\d Z(t) = 2 \sqrt{Z(t)} \d {\gamma}(t),\label{defbessel}\end{equation}
where $\gamma$ is a standard Brownian motion. The absence of drift in the
previous stochastic differential equation makes the function $S_Z$ defined by
$S_Z(x)=x$ for all $x\ge 0$ a scale function of $Z$ (one of
many!). Accordingly, the left-hand side of
(\ref{bess1}) is the probability that, starting from $1$, $Z$ hits $u$ before
hitting 0. This equals $(S_Z(1)-S_Z(0))/(S_Z(u)-S_Z(0)) = 1/u$, proving (\ref{bess1}).\\
 
As for (\ref{bess2}), $Z$ is a martingale whose increasing process is
$\d <Z,Z>_t=4 Z(t) \d t$. Thus, by means of the Dubins-Schwarz theorem
(cf. \cite{R-Y}, p 182), there exists a Brownian
motion, say $\gamma^{*}$, starting from 0, such that for all $t\ge 0$,
$$
Z(t)-1=\gamma^{*}\left( 4\int_0^t Z(s) \d s\right).
$$
\noindent Setting \begin{eqnarray*} {\alpha}_{\delta} &\stackrel{def}{=}& \inf\{s>0:\: 
|Z(s)-1|> \delta\},\\
{\sigma}^{*}_{\delta} &\stackrel{def}{=}& \inf\{s>0:\: |\gamma^{*}(s)|> \delta\},\end{eqnarray*}
\noindent we get that
$${\sigma}^{*}_{\delta}=  4 \int_0^{{\alpha}_{\delta}} Z(s)
  \d s \le 4(1+\delta)\:{\alpha}_{\delta}.$$

\noindent Thus,
\begin{eqnarray*}\p\left(\sup_{0\leq t\leq v} |Z(t)-1| > \delta\right) &=&  
\p\left({\alpha}_{\delta} < v\right) \leq \p\left({\sigma}^ {*}_{\delta}< 4(1+\delta)v\right),\\
&=&\p\left(\sup_{0\leq s \leq 4(1+\delta)v} |\gamma^{*}(s)| >\delta\right),\\
&\leq& 2 \p\left(\sup_{0\leq s \leq 4(1+\delta)v} \gamma^{*}(s)>\delta\right),\\
&\leq& 4 \frac{\sqrt{(1+\delta)v}}{\delta} \exp\left(-\frac{\delta^2}{8(1+\delta) v}\right),\end{eqnarray*} 
as desired. We now move to the 

\bigskip

\noindent {\bf Proof of Lemma \ref{due}:} We start with (\ref{Sigma}). From the definition of
$\Sigma(r)$, it is easily seen that for all $r>0$, 
\[
-\sup_{0\le s \le r} B(s)+ \log\frac{2}{\kappa} (1-e^{-\kappa r/2}) \le \log
 \; (\Sigma(r) e^{-\kappa r /2}) \le -\inf_{0\le s \le r} B(s)+ \log \frac{2}{\kappa}(1-e^{-\kappa r/2}),\]

\noindent and since $\kappa>1$, we have
\[ 
-\sup_{0\le s \le r} B(s) - \log \kappa \le \log (\Sigma(r) e^{-\kappa r/2})\le -\inf_{0\le s \le r} B(s)+ 1\]

\noindent for $r$ large enough. Therefore, for such $r$,  
\begin{eqnarray*}\p(|\log \Sigma (r) -\kappa r /2|>\delta r)&\le& \p(\sup_{0\le s \le
  r}(-B(s))>\delta r/2)+ \p(\sup_{0\le s \le
  r}B(s)>\delta r/2),\\
&=&4 \p(B(1)>\delta {\sqrt r}/2 ) \le 4 \exp(-\frac{\delta^2 }{8} r).
\end{eqnarray*}

\noindent We have used the reflection principle together with a Brownian scaling in deriving the equality above. This finishes the proof of (\ref{Sigma}). 

\bigskip

The next task is to derive (\ref{lam3}) from (\ref{Sigma}). Since $R$ is a Bessel process of dimension $d$ starting at 2, according to Lamperti's
representation (see \ref{Lamperti}), $R$ can be realized as
$$
R(t) = 2 \exp\left(- {1\over 2} \, W_d(\Sigma_d^{-1}(t))\right) , \qquad
t\ge 0, 
$$ 

\noindent where 
$$
W_d(t)= {\gamma}(t) -{d-2\over 4} \, t, \qquad t\ge 0,
$$

\noindent with $\{{ \gamma}(t)\}_{t\ge 0}$ a standard Brownian
motion, and ${\Sigma}_d(t)\; \ddef \; \int_0^t \ee^{-W_d(s)}\d s$. Accordingly, 
\[\int_0^r \frac{\d x}{R^2 (x)} = \frac{1}{4} \Sigma_d^{-1} (r).\]

\noindent Using the above identity, we get that for all $0<\delta<1/(d-2)$, ($d>2$), the left-hand side of 
(\ref{lam3}) equals

\[\p\left(\frac{\log \Sigma_d (s)}{s}-\frac{d-2}{4}<-\frac{\delta(d-2)^2}{8}\right)+\p\left(\frac{\log \Sigma_d (t)}{t}-\frac{d-2}{4}>\frac{\delta(d-2)^2}{4}\right),\]

\noindent where $s=4(\delta+1/(d-2)) \log r$ and $t=4(-\delta+1/(d-2)) \log r$. Making $\kappa=(d-2)/2$ in (\ref{Sigma}), one gets that the probability term in (\ref{lam3}) is less than or equal to $c_3$ times $r^{-c_4}$, where $c_3$ and $c_4$ are constants depending on $d$ and $\delta$. This finishes the proof of (\ref{lam3}) and thus that of Lemma \ref{due}.

\bigskip

\noindent {\bf Proof of Lemma \ref{l:Bessel-ub}:} By scaling, we
can assume without loss of generality that $x=2$. Write for any
$a\ge 0$,
\begin{equation}
    Y_a\; \ddef \; \int_0^\infty {s^a \over R^b(s)} \d s.
    \label{Y_a}
\end{equation}

\noindent We first study the variable $Y_0$. By exactly the same means as in the proof of (\ref{lam3}), namely Lamperti's representation for $R$, one has 

$$
Y_0 = \int_0^\infty {\! \d s\over R^b(s)} = 2^{-b} \int_0^\infty \d {\Sigma}_d (u) \ee^{bW_d(u)/2} = 2^{-b} \int_0^\infty
\ee^{(b-2)W_d(u)/2} \d u .
$$

\noindent By scaling, this implies that 
\begin{equation}
    Y_0 \; \law \; {1\over 2^{b-2} (b-2)^2} \int_0^\infty
    \ee^{{\gamma}(t)-(d-2)t/2(b-2)} \d t .
    \label{Y_0}
\end{equation}

\smallskip

\noindent An application of (\ref{Gamma}) with
$\kappa=(d-2)/(b-2)$ confirms that 
$$
\e_2(Y_0^p)<\infty \; \Longleftrightarrow \; p<{d-2\over b-2}.
$$

\smallskip

\noindent Now consider the variable $Y_a$. For any $t>0$, 
$$
Y_a=\int_0^t {s^a \over R^b(s)} \d s + \int_t^\infty {s^a \over
R^b(s)} \d s \le t^a Y_0 + \int_t^\infty {s^a \over
R^b(s)} \d s.
$$

\noindent For each $p\ge 0$, there exists $d_1(p)$ such that
\begin{equation}
    (x+y)^p \le d_1(p) \, (x^p+y^p), \qquad x\ge 0, \; y\ge 0.
    \label{triangle}
\end{equation}

\noindent Therefore
$$
Y_a^p \le d_1(p) \, t^{ap}Y_0^p + d_1(p) \, \left( \, \int_t^\infty
{s^a \over R^b(s)} \d s\, \right)^p.
$$
\noindent Recall that $\p_x$ and $\e_x$ denote probability and expectation w.r.t. $\p_x$ when the process starts at $x$. 
Taking expectations with respect to $\p_2$ on both
sides and using the Markov property, we obtain
\begin{equation}
    \e_2(Y_a^p) \le d_1(p) \, t^{ap} \, \e_2(Y_0^p) + d_1(p) \,
    \e_2 \left( \, \e_{R(t)} \left( \, \int_0^\infty {(s+t)^a \over
    R^b(s)} \d s\, \right)^p \, \right).
    \label{Y_a^p}
\end{equation}

\medskip

\noindent According to (\ref{triangle}), 
\begin{eqnarray*}
    \left( \, \int_0^\infty {(s+t)^a \over R^b(s)} \d s\, \right)^p
 &\le& d_1^p(a)\, \left( \, \int_0^\infty {s^a \over R^b(s)} \d s +
    t^a \int_0^\infty {\! \d s\over R^b(s)} \right)^p
\cr &\le& d_2 \, \left( \left( \, \int_0^\infty {s^a
    \over R^b(s)} \d s\right) ^p  + t^{ap} \left( \, \int_0^\infty
    {\! \d s\over R^b(s)} \right)^p \right),
\end{eqnarray*}

\noindent where $d_2=d_2(a,p)\; \ddef \; d_1^p(a) d_1(p)$. Applying the
scaling property yields that for any $y>0$, 
$$
\e_y \left( \, \int_0^\infty {(s+t)^a \over R^b(s)} \d s\,
\right)^p \le d_2 \, \left( (y/2)^{(2a-b+2)p} \e_2 (Y_a^p) +
(2/y)^{(b-2)p} t^{ap} \, \e_2(Y_0^p) \right).
$$

\smallskip

\noindent Plugging this into (\ref{Y_a^p}) gives that for any
$t>0$,  
\begin{eqnarray*}
    \e_2(Y_a^p) 
 &\le& d_1(p) \, t^{ap} \, \e_2(Y_0^p) + d_3 \, \e_2(
    R^{-(b-2a-2)p}(t))\, \e_2 (Y_a^p)
\cr && \qquad + d_4 \, t^{ap} \, \e_2( R^{-(b-2)p}(t))\,
    \e_2 (Y_0^p) ,
\end{eqnarray*}

\smallskip

\noindent with $d_3=d_3(a,b,p)\; \ddef \; 2^{-(2a-b+2)p} d_2$ and
$d_4=d_4(a,b,p)\; \ddef \; 2^{(b-2)p} d_2$. For any $0<u<d$,
$$
\e_2( R^{-u}(t)) \le {d_5(u,d)\over t^{u/2}}, \qquad t\ge 1,
$$

\noindent for some $d_5(u,d)$ (this can be easily checked for
example using the exact semi-group of $R$). Therefore, if
$0<p<(d-2)/(b-2)$ which guarantees $0<(b-2a-2)p<d$, then we can 
choose $t$ sufficiently large that $\e_2( R^{-(b-2a-2)p}(t))
\le 1/(2d_3)$, which yields
$$
\e_2(Y_a^p) \le 2d_1(p) \, t^{ap} \, \e_2(Y_0^p) + 2 d_4 \,
t^{ap} \, \e_2( R^{-(b-2)p}(t))\, \e_2 (Y_0^p) .
$$

\noindent In particular, this shows that $\e_2(Y_a^p)< \infty$ for
all $p<(d-2)/(b-2)$. 

\bigskip

On the road to the proof of Theorem \ref{t:main}, our first step is to prove Theorem \ref{t:H}
for the first hitting time process $H$. This will be justified in Section
\ref{s:main}. The following section provides our second step, as finding tail estimates for $H$ amounts to
finding tail estimates for two random variables $I_1$ and $\Upsilon$, to be defined below.

\section{From hitting times to Bessel and Jacobi processes}
\label{s:Second reduction}
Recall the definitions of $\{H(r)\}_{r>0}$ and $\{\sigma_{{\cal
    B}} (r)\}_{r>0}$ from (\ref{H}) and (\ref{sigma}), where ${\cal B}$ is a
    Brownian motion independent of the environment $W$, see (\ref{Brox}). In
    the sequel, we shall drop ${\cal B}$ in both $\sigma_{{\cal
    B}} (r)$ and $L^x_t({\cal B})$ for brevity. By (\ref{Brox}) and the occupation density formula we have, for any $r>0$,
\begin{eqnarray*}
H(r)= T(\sigma (S(r))) &=& \int_0^{\sigma(S(r))}
\ee^{-2W(S^{-1}({\cal B}(u)))} \d u, \\
&=& \int_{-\infty}^{S(r)}  e^{-2W(S^{-1} (y))} L_{\sigma (S(r))}^y \d y,\\
&=& \left(\int_{-\infty}^0 + \int_0^r\right) e^{-W(x)} L_{\sigma
(S(r))}^{S(x)} \d x, \\
&\stackrel{def}{=}& I_1 (r)+I_2 (r),
\end{eqnarray*}

\noindent where we have performed the change of variables $x=S^{-1} (y)$ in deriving the third equality.

\medskip

A scaling argument tells us
that, at fixed environment $W$, the processes $\{L_{\sigma(S(r))}^{S(y)}\}_{y
  \le r}$ and $\{S(r) 
\: L_{\sigma(1)}^{S(y)/S(r)} \}_{y\le r}$ have the same law.

\smallskip
Further, according to the first 
Ray--Knight theorem (see Fact 1, Section \ref{s:preliminaries}), $\{L_{\sigma(1)}^{-z}\}_{z\ge 0}$ is a squared Bessel process
of dimension $0$, with initial exponential distribution $\xi$ of mean $2$ and $\{L_{\sigma(1)}^{1-t}\}_{0\leq t \leq 1}$ 
is a two-dimensional squared Bessel process starting from 0 and independent of $W$, say $R_1^2$. 

\medskip

As result, $H(r)$ can be rewritten as:
\begin{eqnarray}\label{I_1-law}
H(r)&\stackrel{law}{=}& S(r) \, \xi \int_{-\infty}^0 \ee^{-W(y)} \, Z\left( { |S(y)| \over S(r) \, \xi}\right) \d y 
+  S(r) \, \int_0^r R_1^2 \left(1-\frac{S(y)}{S(r)}\right) \d y,
\end{eqnarray}

\noindent where $\{Z(t)\}_{t\ge 0}$ is a squared Bessel process of dimension $0$, starting from 1, and where $S(r)$, $\{ W(y)\}_{y\le 0}$, $\xi$ and $\{Z(t)\}_{t\ge 0}$ are independent. 

\medskip

Note that the above identity in law is {\it quenched}, hence also annealed. 

\medskip

A glance at the definition of $I_1 (r)$ tells us that $r
\mapsto I_1 (r)$ increases so that, $\p$-almost surely, 

\begin{equation}\label{I_1fini}\sup_{r\ge0} I_1 (r)=I_1 (\infty)<\infty\end{equation}
Actually, since $\kappa>1$, both $S(\infty)$ and $\int_{-\infty}^0
e^{-W(y)} \d y$ have finite expectations; it is then easily checked, from (\ref{I_1-law}), that $\e(I_1 (\infty))<\infty$. 

\bigskip

Now, in order to estimate the tail probabilities of $I_2 (r)$, a slight
transformation of the expression given in (\ref{I_1-law}) is needed. 

\medskip

At fixed $W$, a scaling argument used twice, followed by the change of variables $z=r-y$, leads to the following series of quenched identities in law:

\begin{eqnarray*} I_2 (r) &\stackrel{law}{=} & S(r) \; \int_0^r e^{-W(y)} R_1^2 \left(1-\frac{S(y)}{S(r)}\right) \; \d y,\\
&\stackrel{law}{=} & \int_0^r e^{-W(y)} R_1^2 (S(r)-S(y))\; \d y,\\
&\stackrel{law}{=} &  \int_0^r e^{W(r)-W(y)} \;  R_1^2 \left(\int_y^r e^{-(W(r)-W(x))} \; \d x\right),\\
&=& \int_0^r e^{B^r_z-\kappa z/2} R_1^2 \left(\int_0^z  e^{-(B^r_x-\kappa x/2)} \d x\right)\; \d z,
\end{eqnarray*} 

\smallskip
\noindent where $B^r_z =B(r)-B(r-z)$, for $0\le z \le r$. Since $B^r$ and $B$ have the same law on $[0,r]$, one gets the following 
{\it annealed} identity in law:
\begin{eqnarray*} I_2 (r)&\stackrel{law}{=} & \int_0^r e^{W(y)} R_1^2\left(\int_0^y e^{-W(x)} \d x\right) \; \d y.\end{eqnarray*}

\noindent Using Lamperti's representation (see (\ref{Lamperti}) with $\zeta=\kappa$) gives that

\[
e^{-W(y)}= {1\over 4} \, R_2^2 \left(\int_0^y e^{-W(x)} \d x\right),\]

\noindent where $R_2$ is a transient Bessel process of dimension $2+2\kappa$, starting from 2. In this light, denoting for simplicity $\int_0^y e^{-W(x)} \d x$ by $\Sigma(y)$, then performing the change of variables $u=\Sigma(y)$, we arrive at:

\begin{equation}\label{defI3}I_2 (r)\stackrel{law}{=}
16\int_0^{\Sigma(r)}\frac{R_1^2(u)}{R_2^{4}(u)}\d u\stackrel{def}{=}16I_3 (\Sigma(r)).\end{equation}

\noindent Observe that $R_2$ depends only on the environment $W$ (of which $R_1$ is independent).

\bigskip

According to a result by Warren and Yor, see (\ref{warren}), there exists a Jacobi process of dimensions $(2,2+2\kappa)$, say $Y$, starting from 0, such that 
\begin{equation}\label{I_3}
I_3 (r)\stackrel{law}{=}\int_0^{\Lambda(r)} \frac{Y(s)}{(1-Y(s))^2} \d s\; \ddef\; \Upsilon(\Lambda(r)),
\end{equation}
\noindent where 
\begin{eqnarray}\label{upsilon}
\Upsilon(r)\; \ddef\; \int_0^r \frac{Y(s)}{(1-Y(s))^2} \d s, \\
\Lambda(r)\; \ddef\; \int_0^r \frac{\d u}{R_1^2(u)+R_2^2(u)}.\label{Lambda}\end{eqnarray}

\indent Note that since $R^2\stackrel{def}{=}R_1^2+R_2^2$, the process $R$ is a squared Bessel process of 
dimension $d\stackrel{def}{=}4+2\kappa >2$, starting from $2$. We know from (\ref{lam1}) that 
$\Lambda(r)/ \log r$ approaches, $\p$-almost
surely, $1/(d-2)= 1/(2+2\kappa)$. Moreover, Lemma \ref{due} makes us expect that tail estimates for $I_2$ will follow from those for 
$\Upsilon$. And they do:

\begin{theorem}\label{p:derniere}
\noindent For all $u>(1+\kappa)/(\kappa (\kappa-1))$,
\begin{equation}\label{upsilon:ub}\limsup_{r\rightarrow \infty} \frac{1}{\log r} \log \p\left(\Upsilon(r)>u r\right) \le 1-\kappa,\end{equation}
\noindent and for all $u>0$, 
\begin{equation}
\label{upsilon:lb}\liminf_{r\rightarrow \infty} \frac{1}{\log r} 
\log \p\left(\Upsilon(r)>u r\right) \ge 1-\kappa.\end{equation}

\noindent Moreover, for any $w>0$,
   \begin{equation}
      \limsup_{r\to \infty} \frac{1}{\log r} \log \p( I_1(\infty)>wr)
      \leq 1-\kappa.
      \label{I_1-tail}
   \end{equation}
\end{theorem}  

\bigskip

\noindent {\bf Proving that Theorem \ref{t:H} follows from Theorem \ref{p:derniere}:} Assuming that both (\ref{upsilon:ub}) and (\ref{I_1-tail}) hold true, we first
prove (\ref{t:H-ub}), the upper bound for $H$. 

\medskip

\noindent To this end, let $u>v_{\kappa}^{-1}=4/(\kappa-1)$ be given. Then there exists some $\delta_0$ such that 
$u>(1+\delta_0)v_{\kappa}^{-1}.$ Next pick $\epsilon$ with $0<\epsilon<\delta_0 / (1+\delta_0)$ so that $(1-\epsilon)(1+\delta_0) > 1$, then choose $\delta>0$ so small that $(1-\epsilon) (1+\delta_0) > 1+2\delta$. 

\medskip

Now, since $H(r)=I_1 (r)+I_2(r)$, using successively the triangle inequality,
(\ref{I_1fini}), the definition of $I_3$, (\ref{defI3}), the fact that
$r\mapsto I_3 (r)$ is increasing, then finally Lemma \ref{due}, 
(\ref{Sigma}), it follows that 
\begin{eqnarray*}
\p(H(r)>ur)&\le& \p(I_1 (r)>\epsilon u r)+\p(I_2(r) >(1-\epsilon) ur),\\
&\leq& \p(I_1 (\infty)>\epsilon u r)+\p(I_2(r) >(1-\epsilon) ur),
\end{eqnarray*}

\noindent with
\begin{eqnarray*}\p(I_2(r) >(1-\epsilon) ur)&=& \p\left(16 \; I_3(\Sigma(r))>(1-\epsilon)ur\right),\\
&\le&\p\left(I_3(e^{\kappa(1+2\delta)r/2})>\frac{(1-\epsilon)}{16} ur\right) +c_1e^{-c_2 r},\\
&=&\p\left(I_3(s)>v\log s\right)  +c_1 e^{-c_2 r},
\end{eqnarray*}

\smallskip
\noindent where $s=s(r)=e^{\kappa(1+2\delta)r/2}$ and 
$$
v=\frac{1-\epsilon}{1+2\delta} \times \frac{u}{8\kappa}> \frac{1}{2
  \kappa(\kappa-1)} \times \frac{(1-\epsilon)(1+\delta_0)}{1+2\delta} > \frac{1}{2\kappa(\kappa-1)},$$
\noindent thanks to the choices of $\delta_0$, $\epsilon$ and $\delta$. 

\medskip

Now, since $v>(2\kappa (\kappa-1))^{-1}$, there exists $0<\epsilon_0 <1$ such
that $v>(1+\epsilon_0)/2\kappa (\kappa-1)$. Since $4+2\kappa >2$, (\ref{lam3}) 
tells us about the rate at which $\Lambda(s)/\log s$ approaches $\p$-almost
surely $1/(2+2\kappa)$ as $s$ tends to infinity. Knowing (\ref{I_3}), this in conjunction with the fact that $\Upsilon$ is increasing yields
\begin{eqnarray*}
\p(I_3(s)>v \log s)&\le& \p(\Upsilon(t)>v\log s)+\p\left
(\left|\frac{\Lambda(s)}{\log s}-\frac{1}{2+2\kappa}\right|> \frac{\epsilon_0}{2+2\kappa} \right),\\
&\le& \p(\Upsilon(t)>wt)+\frac{c_3}{s^{c_4}},
\end{eqnarray*}

\noindent for all $n\ge 1$, and for some constants $c_3$ and $c_4$ depending on $\epsilon_0$ and $\kappa$, where 
$$t=t(s)=\frac{1+\epsilon_0}{2+2\kappa} \log s, \; \; \; \; 
w=\frac{2+2\kappa}{1+\epsilon_0} v>\frac{1+\kappa}{\kappa(\kappa-1)}.$$

\noindent Consequently, putting all the pieces together, and keeping in mind that $t$ is actually $r$ times a 
constant depending on $\delta$, $\epsilon_0$ and $\kappa$, one gets
\begin{equation}\label{allez}
\p(H(r)>ur)\le \p(I_1(\infty)>\epsilon u r)+\p(\Upsilon(t)>wt)+ c_3 e^{-c_4 r/2}
,\end{equation}

\noindent for $\kappa>1>\delta>0$. Taking the logarithm of both sides of (\ref{allez}), using
the elementary fact that $\log(a+b+c)\le \log 3+\sup(\log a, \log b, \log c)$
for $a,b,c>0$, dividing by $\log r$ ($\log r \sim \log t$, as $r\rightarrow
\infty$), taking the limsup, and making use of (\ref{upsilon:ub}) and
(\ref{I_1-tail}) completes the proof of the upper bound for $H$ (\ref{t:H-ub}). 

\medskip

As for the lower bound for $H$ (\ref{t:H-lb}), this follows from the lower
bound for $\Upsilon$ (\ref{upsilon:lb}), by Lemma \ref{due}. The reasoning is the same as before (but a bit simpler actually since $I_1 (\infty)$ does not enter the picture): we write $\p(H(r)>ur)\geq \p(I_2 (r)>ur)=\p(16 I_3 (\Sigma(r))>ur),$ and use the same arguments as before.

%we give just a sketch of the proof. Using the same notations as before, for any $u>0$,
%we have
%\begin{eqnarray*}\p(H(r)>ur)&\geq &\p(I_2 (r)>ur)=\p(16 I_3 (\Sigma(r))>ur),\\
%& \geq& \p(I_3 (s)>v^{'} \log s)-c_1 e^{-c_2 r},\\
%&\geq & \p(\Upsilon(t)>w^{'} t)-c e^{-c^{'} r},
%\end{eqnarray*}
%\noindent thanks to the same arguments as for the upper bound. 

\smallskip

This indicates how (\ref{t:H-lb}) follows from (\ref{upsilon:lb}).

\bigskip

\noindent We have seen that estimating tail probabilities for $H(r)$ reduces to proving Theorem \ref{p:derniere}. 
This amounts to studying the tail asymptotics for $I_1(\infty)$ and $\Upsilon$. We postpone the study of $I_1(\infty)$ 
to Section \ref{s:End} and move on to the proofs of (\ref{upsilon:ub}) and (\ref{upsilon:lb}) of Theorem \ref{p:derniere}. In the next section, we state and prove a key result which will
enable us to prove these results. 

\bigskip

\section{A key estimate}\label{s:key}
From the stochastic differential equation (\ref{jacobieds}), we have
%a glance at \cite{K-T} pages 194-195  gives:
\[\d S_Y (x)= \frac{\d x}{x(1-x)^{\kappa+1}}, \; \; \; 0<x<1,\]
\noindent and
\[m_Y (dx)=\frac{1}{4} (1-x)^{\kappa} \d x,\; \; \; 0<x<1,\]

\medskip
\noindent where $S_Y$ is a scale function of the diffusion $Y$ and $m_Y$ is its speed measure. 

\medskip
 \noindent Recall that $Y_0 =0$ and, given the definitions of $R_1$ and $R_2$, that $0<Y(t)<1$ for all $t>0$, by (\ref{warren}).

\medskip
 It is easily checked that $Y$ is recurrent, and that 0 and 1 are actually entrance boundaries 
(see \cite{K-T}, page 235); they cannot be reached from ]0,1[. Since in our
 case $Y$ starts at 0, it {\it rapidly} moves to ]0,1[ never to return to 0. Tail estimates for the first hitting time of level 1/2 by the diffusion $Y$, started at 0, 
shed some light on just how fast $Y$ moves from 0 to 1/2; see Lemma (\ref{l:quick}). 

\medskip

In order to establish (\ref{upsilon:ub}) and (\ref{upsilon:lb}) of Theorem \ref{p:derniere}, we first assume that $Y$ starts in ]0,1[, choosing without loss of generality that it starts at 1/2, and get the desired estimates with $\p_{1/2}$ replacing $\p$. Next, as is proved in the next section, Lemma \ref{l:quick} in conjunction with the strong Markov property enables us to establish the result for the case where $Y$ 
starts at 0; this then yields (\ref{upsilon:ub}) and (\ref{upsilon:lb}). 

\bigskip 

Since $Y(0)=1/2$, a scale function of $Y$ is
\begin{equation}\label{SY}S_Y(x)=\int_{1/2}^y\;  \frac{dx}{x(1-x)^{\kappa+1}}, 
\; \;  y\in(0,1),\end{equation}
so that \begin{equation}\label{S(Y)}dS_Y(Y(t))=2 Y^{-1/2}(t) (1-Y(t))^{-\kappa-1/2} 
d\b(t).\end{equation}

Thus, $Y$ can be constructed from a Brownian motion via a scale transformation and time change. Namely, there exists a driftless Brownian motion $\beta$ such that
\begin{equation}S_Y(Y(t))=\beta(U(t)) \label{itoMc}, \; \; \; t\geq 0,\end{equation}

\noindent where the time change $U$ is given by 
\begin{equation}\label{UU}U(t)\; \ddef\;  4\int_0^t \frac{ds}{Y(s)(1-Y(s))^{2\kappa+1}}.\end{equation}

\bigskip

From \cite{HSY} we get that, for a certain Brownian motion $\beta_r$ defined
below, see (\ref{beta}), $U(r)/r^2$ is roughly $\tau_{\beta_r}(4(\kappa+1))$.
The next proposition provides a key estimate which measures the
error introduced by this replacement. Before stating this, we 
need some definitions. Let
$\alpha,\; \delta,\; \mu,\; \nu>0$ such that $\delta<1$ and $\mu + \nu <1$, and set  
\begin{equation}\label{lambda+}
\lambda_{\pm} \; \ddef\; \lambda(1\pm \delta), \qquad {\rm where} \quad \lambda\; \ddef\;  4(\kappa+1),
\end{equation}

\begin{equation} \label{choices}0<\theta<1, \; \;  \mu=e^{-r^{\theta}}\; \;  {\mathrm and}\; \;  \nu = 
r^{-\theta/{\kappa}},\end{equation}

\begin{equation}\label{beta} \beta_r (s)\; \ddef\; \frac{1}{r}\beta(r^2 s), \qquad s\ge 0,\end{equation}

\begin{equation}\label{epsilon}\epsilon(r,s)\; \ddef\;\frac{1}{4} \int_0^1 (1-x)^{\kappa}\; \left(L_s^{S_Y(x)/r}(\beta_r)-L_s^0(\beta_r)\right) \d x.\end{equation}

\noindent and 
\begin{eqnarray}\label{Ki}\Xi_{r}=\Xi(r, \kappa,\delta)\;\ddef\; \left\{
\epsilon(r,\tau_{\beta_r}(\lambda_+)) \ge -\delta\; ; \; 
\epsilon(r,\tau_{\beta_r}(\lambda_-))\le \delta\right\}. \end{eqnarray}

\bigskip

So as not to overburden the reader with notation, we shall drop $\beta_r$ in both $\tau_{\beta_r} (.)$ and $L_{\tau_{\beta_r}
  (.)}^. (\beta_r)$ throughout.

\bigskip

\begin{Proposition}\label{star}
Let $0<\delta<1$ be given. On the event $\Xi_{r}$, we have, for all $r>0$, 
\begin{equation}\label{ineg}
{\tau}(\lambda_-)\leq\frac{U(r)}{r^2}\leq
{\tau}(\lambda_+),
\end{equation}

\noindent and 
\begin{equation}\label{probastar}\lim_{r \rightarrow \infty} \frac{\log
    \p_{1/2}\left(\Xi_{r}^c\right)}{\log r}=-\infty,\end{equation}

\smallskip

\noindent where $\Xi_{r}^c$ stands for the complement of $\Xi_{r}$.
\end{Proposition}

\bigskip

\noindent {\bf Proof of Proposition \ref{star}} 
Let $U^{-1}$ denote the inverse of $U$. Thanks to (\ref{itoMc}), the
density occupation formula (\ref{occupation}), and (\ref{UU}), we arrive (exactly as in \cite{HSY} with $d_1=2$ and $d_2=2+2\kappa$) at the following:
\begin{eqnarray}U^{-1}(t)&=&\frac{1}{4} \int_0^t \left(S_Y^{-1}(\beta(s))\right)\; 
\left(1-S_Y^{-1}(\beta(s))\right)^{1+2\kappa} \d s \nonumber \\
&=& \frac{1}{4} \int_0^1 (1-x)^{\kappa} L_t^{S_Y(x)}(\beta) \d x.\nonumber 
\end{eqnarray} 
Hence, 
\begin{eqnarray*}
\frac{1}{r} U^{-1}(r^2 s)&=&\frac{1}{4} \int_0^1 (1-x)^{\kappa}\; L_s^{S_Y(x)/r}(\beta_r) \d x,\\
&=& \frac{1}{\lambda} L_s^0 (\beta_r)+\epsilon(r,s).\end{eqnarray*}

\noindent Since $U$ is increasing, (\ref{Ki}) together with
straightforward computations delivers (\ref{ineg}). 

\bigskip

We now turn to (\ref{probastar}). Recalling from (\ref{lambda+}) and
(\ref{epsilon}) the definitions of $\lambda_{+}$ and $\epsilon(.,.)$ respectively, we have that for all $\mu$ and $\nu$ positive such that $\mu
+ \nu <1$, 
\begin{eqnarray*}
\p_{1/2}\left( \epsilon(r, \tau (\lambda_{+}))<-\delta\right)
&\le& \p\left(\int_{\mu}^{1-\nu} (1-x)^{\kappa}
  L_{\tau(\lambda_{+})}^{S_Y(x)/r} \d x<4\right),\\ 
&\le& \p\left(\inf_{\frac{S_Y (\mu)}{r \lambda_{+}} \le x \le \frac{S_Y(1-\nu)}{r \lambda_{+}}} L_{\tau(1)}^x <1-\delta_1\right),
\end{eqnarray*}

\noindent where 
$$
\delta_1=\delta_1(\delta, \mu, \nu, \kappa)=1-\frac{1}
{(1+\delta)\left((1-\mu)^{\kappa+1}-{\nu}^{\kappa+1}\right)}.$$ 

\noindent In deriving the last inequality, we have used a scaling argument together with the monotonicity of
$S_Y$. A little Brownian excursion
theory now tells us that $\{L_{\tau(1)}^x\}_{x \ge 0}$ and
$\{L_{\tau(1)}^{-x}\}_{x \ge 0}$ are two independent squared Bessel processes
of dimension 0, started at 1, (see \cite{R-Y}). Therefore the last probability above is  
\begin{eqnarray*}
&\le&2 \: \p\left(\inf_{0 \le x \le \frac{S_Y(1-\nu)\vee 
|S_Y (\mu)|}{r \lambda_{+}}}
  L_{\tau(1)}^x < 1-\delta_1\right), \\&\le&
2 \: \p\left(\sup_{0 \le x\le \frac{S_Y(1-\nu)\vee |S_Y (\mu)|}{r\lambda_{+}}} |L_{\tau(1)}^x -1|>\delta_1\right).\end{eqnarray*}

\noindent Recalling the definition of $S_Y$, (\ref{SY}), straightforward computations tell us that
\begin{eqnarray} 2^{\kappa +1} \log 2x &\leq& S_Y (x) \leq \frac{1}{(1-x)^{\kappa+1}} \: \log
2x, \nonumber\\
\label{truc}\frac{1}{ 2 \kappa x^{\kappa}} &\leq& S_Y (1-x)\leq \frac{2}{x^{\kappa}},\end{eqnarray}

\noindent as $x$ approaches 0.

\medskip

\noindent Thus, with the choices of $\mu$ and $\nu$, see (\ref{choices}), we get 
$$\frac{1}{r \lambda_{+}}(S_Y(1-\nu)\vee |S_Y (\mu)|)\le 2^{\kappa-2} r^{\theta-1},$$

\noindent for $r$ large enough. Accordingly, with the help of (\ref{bess2}),
%in which case $\delta_1$ is in the vicinity of
%$\delta/(1+\delta)$, and hence $0<\delta_1<1$. 
\begin{eqnarray*}  2\p\left(\sup_{0 \le x\le \frac{S_Y(1-\nu)\vee 
|S(\mu)|}{r
      \lambda_{+}}}   |L_{\tau(1)}^x -1|>\delta_1\right)& \le & 2\p\left(\sup_
{0 \le x\le 2^{\kappa-2} r^{\theta-1}}
    |L_{\tau(1)}^x -1|>\delta_1\right), \\
&\le& f_1 \ee^{- f_2 r^{1-\theta}},\end{eqnarray*}

\noindent where $f_1= 5\: 2^{\kappa/2}/\delta_1$ and $f_2=\delta_1^2 \:
      2^{-\kappa-2}$. 

\medskip

\noindent As a result,

\begin{equation}
\lim_{r \rightarrow \infty} \frac{1}{\log r} \log \p_{1/2}\left(\epsilon(r,
\tau({\lambda}_{+}))<-\delta\right)=-\infty.\label{mimi0}\end{equation}

\bigskip
\medskip

We now turn to $\p_{1/2}\left(\epsilon(r,\tau({\lambda}_{-}))>\delta\right)$. It is plain to see that
 %That $\lambda_{-}$
%plays no role in the proof is obvious. Indeed, thanks to a scaling argument,
%and upon replacing $r$ and $\delta$ by $r\lambda_{-}$ and
%$\delta/{\lambda_{-}}$ respectively, all we need to upperbound is $\p_{1/2}\left(\eps%ilon(r,\tau(1),\kappa)>\delta\right).$
\begin{equation}
\p_{1/2}\left(\epsilon(r,\tau(\lambda_{-})) > \delta\right)\le\p_{1/2}(J_1(r, \mu)> 2\delta)+ \p_{1/2}(J_2(r, \mu, \nu)> \delta)+\p_{1/2}(J_3(r, \nu)> \delta).
\end{equation}

\noindent where\begin{eqnarray*}J_1(r,\mu)&\stackrel{def}{=}&\int_0^{\mu} 
(1-x)^{\kappa}\left(L_{\tau (\lambda_{-})}^{S_Y(x)/r}-\lambda_{-}\right)\d
x,\nonumber \\
J_2(r,\mu, \nu)&\stackrel{def}{=}&
\int_{\mu}^{1-\nu} (1-x)^{\kappa}\left(L_{\tau (\lambda_{-})}^
{S_Y(x)/r}-\lambda_{-} \right)\d x,\nonumber \\
J_3(r,\nu)&\stackrel{def}{=}&\int_{1-\nu}^1 (1-x)^
{\kappa}\left(L_{\tau (\lambda_{-})}^{S_Y(x)/r}-\lambda_{-} \right)\d x.\end{eqnarray*}

\noindent We start with $\p_{1/2}(J_1(r, \mu)> 2\delta)$. Since $\sup_{x\le
  0}L_{\tau(1)}^x$ and $\sup_{x\ge
  0}L_{\tau(1)}^x$ have the same law, (\ref{bess1}) together with a scaling
  leads to
\begin{eqnarray}
\p_{1/2}(J_1(r, \mu)> 2\delta)&\le& \p\left(\sup_{x\le 0}L_{\tau(1)}^x
> \frac{\delta}{\lambda_{-}\mu}\right),\nonumber \\
&=&  \p\left(\sup_{x\ge 0}L_{\tau(1)}^x > \frac{\delta}
{\lambda_{-} \mu}\right) = \frac{\lambda_{-}}{\delta} \mu=\frac{\lambda_{-}}{\delta} \: e^{-r^{\theta}}.\label{mimi1}
\end{eqnarray}

\noindent Next, we may write  
$$\p_{1/2}(J_2(r, \mu, \nu) > \delta) \le 2 \: \p\left(\sup_{0\le x \le (S_Y(1-\nu)\vee |S_Y(\mu)|)/ \lambda_{-}r}\left|L_{\tau(1)}^x-1\right|> \delta_2\right),$$

\noindent where $$\delta_2=\delta_2(\delta, \mu, \nu, \kappa)=\frac{\delta}
{4(1-\delta)\left((1-\mu)^{\kappa +1}- {\nu}^{\kappa +1}\right)}.$$

\smallskip

\noindent Thus, for $r$ large, with the same choice of $\mu$ and $\nu$ as
before, one gets 
%$\delta_2$
%approaches $\delta(\kappa+1)<1$, thus $0<\delta_2 <1$ leading to
\begin{equation}\label{mimi2}
\p_{1/2}(J_2(r, \mu, \nu)> \delta)\leq f_3 \ee^{-f_4 \: r^{1-\theta}}.\end{equation}

\noindent where $f_3= 2^{2+\kappa/2}  \:
\sqrt{1+\delta_2}\delta_2^{-1},$ and
$f_4=\delta_2^2 \: 2^{-\kappa-1}(1+\delta_2)^{-1}$.

\bigskip

\noindent Lastly, we find an upper bound for $\p_{1/2}(J_3(r, \nu)>\delta)$. 
Thanks to (\ref{truc}), $1-S_Y^{-1} (y) \leq (2/y)^{1/\kappa}$ for $y$ large
enough. As a result, performing the change of variables $y=S_Y (x)/r$ together with a
change of scale yields
\begin{eqnarray}
\p_{1/2}(J_3(r,\nu) > \delta)
&\le& \p\left(\int_{1-\nu}^1 \: (1-x)^{\kappa} \: L_{\tau(\lambda_{-})}^{S_Y
    (x)/r} \d x>\delta \right),\nonumber \\
& \le& \p\left(\int_{fr^{\theta-1}}^{\infty} \frac{L_{\tau(1)}^y}{y^{2+1/\kappa}}
 \d y > \delta_3 r^{1+1/\kappa}\right),\label{iteration}
\end{eqnarray}

\noindent where $f=(2\kappa \: \lambda_{-})^{-1}$ and $\delta_3 = 
\frac{\delta}{4} (\lambda_{-}/2)^{1/\kappa}$. 

\medskip

\noindent Now, the definition of a Bessel process of dimension
0, see (\ref{defbessel}), together with the integration by parts
$$
\frac{L_{\tau(1)}^y}{y^{2+1/\kappa}} d y = -d\left( L_{\tau(1)}^y 
\frac{1}
{(1+1/\kappa) y^{1+1/\kappa}}\right)+ 2 \frac{\sqrt{L_{\tau(1)}^y}}
{(1+1/\kappa) y^{1+1/\kappa}} \d {\gamma(y)},$$

\medskip

\noindent implies that the last probability above is
\begin{eqnarray*}
&\le& \p\left(L_{\tau(1)}^{fr^{\theta-1}}> \delta_4
    r^{\theta(1+1/\kappa)}\right)+ \p\left(\int_{fr^{\theta-1}}^{\infty} \frac{\sqrt{L_{\tau(1)}^y}}{y^{1+1/{\kappa}}} \d{\gamma(y)} > \delta_5 r^{1+1/\kappa}\right),
\end{eqnarray*}

\noindent where 
$\delta_4 = \frac{\delta_3}{2} (1+\frac{1}{\kappa}) f^{1+1/\kappa}$ and 
$\delta_5=  \delta_3 (1+1/\kappa)/4$. Recall that $L_{\tau(1)}^y$ almost surely goes to $0$ as $y$ goes to infinity.

\medskip

\noindent An exponential inequality together with (\ref{bass-griffin}) for
$z=fr^{\theta-1}$ and $u=1/(3z)$ tells us that 
$$
\p \left(L_{\tau(1)}^{fr^{\theta-1}}> \delta_4
    r^{\theta(1+1/\kappa)}\right)\le e^{-\frac{\delta_4}{4f} \: r^{1+\theta/\kappa}},$$

\noindent for $r$ large enough. On the other hand, writing the stochastic
integral above as a time changed Brownian motion gives, again for $r$
sufficiently large, 
$$
\p\left(\int_{fr^{\theta-1}}^{\infty}
 \frac{\sqrt{L_{\tau(1)}^y}}{y^{1+1/\kappa}} 
\d {\gamma(y)} \ge  \delta_5 r^{1+1/\kappa}\right)\le e^{-\frac{\delta_5^2}{2}\log^2 r} + 
\p\left(\int_{fr^{\theta-1}}^{\infty} 
\frac{L_{\tau(1)}^y}{y^{2+2/\kappa}}\: \d y> \frac{r^{2+2/\kappa}}{\log^2 r}\right).$$

\noindent The last term is nothing but $\p(\sup_{0\leq y \leq 1} \gamma (y)
> \delta_5 \: \log r)=2\;\p({\cal {N}}> \delta_5 \: \log r)$, thanks to the
  reflection principle, see \cite{R-Y}, for ${\cal {N}}$ a normalized Gaussian variable.

\medskip

\noindent {\bf The iteration scheme:}
\smallskip

\noindent We iterate the procedure above $m$ times, which gives that 
$\p_{1/2}(J_3(r,\nu) > \delta)$ is 
\begin{equation}\label{mimi3}
\le m e^{-f_m \log^2 r}+ m e^{-\frac{\delta_4}{4f} \: r^{1+\theta/\kappa}}+
\p\left(\int_{fr^{\theta-1}}^{\infty}\frac{L_{\tau(1)}^y}{y^{2+2^m
      /{\kappa}}} \d y >  \frac{r^{2^m +2^m /\kappa}}{(\log r)^{2^{m+1}-2} }\right),
\end{equation}

\medskip

\noindent where $f_m$ is a constant depending only on $f$ and the integer $m$, or equivalently,
on $\kappa$, $\delta$ and $m$.

\bigskip

\noindent Since $\{L_{\tau (1)}^y -1\}_{y\geq 0}$ is a
martingale, see for instance \cite{R-Y}, $\e(L_{\tau (1)}^y
)=1$ for all $y\geq 0$, and thus, by Chebychev's inequality  
\begin{equation}\label{mimi4}
\p\left(\int_{fr^{\theta-1}}^{\infty}\frac{L_{\tau(1)}^y}{y^{2+2^m /{\kappa}}} \d y >  \frac{r^{2^m +2^m /\kappa}}{(\log r)^{2^{m+1}-2}}\right)\le 
f^{-1-2^m/\kappa}\; \frac{(\log r)^{2^{m+1}-2}}{r^{2^m (1+\theta/\kappa)+\theta-1}}.\end{equation}

\medskip

\noindent Recall that (\ref{mimi0}) takes care of $\p_{1/2}\left(\epsilon(r,
  \tau (\lambda_{+}))<-\delta\right)$. On the other hand, putting (\ref{mimi1}), (\ref{mimi2}) and (\ref{mimi4}) together gives that 
$$\limsup_{r \rightarrow \infty} \frac{1}{\log r}\log \p_{1/2}\left(\epsilon(r, 
\tau (\lambda_{-}))>\delta\right)
 \le  -\left(2^m(1+\theta/\kappa)+\theta-1\right),$$
for all $0<\theta<1$ and any fixed but arbitrary integer $m$. Letting $m$ go to
 infinity gives that the limsup above is in fact a limit, which equals
 $-\infty$. The proof of (\ref{probastar}) is now complete. We are ready for the 

\medskip

\section{Tail estimates for $\Upsilon$}
\label{s:six}
As announced in the beginning of last section, getting tail estimates for $\Upsilon$ will split into two parts: we first 
deliver the desired result assuming that $Y(0)=1/2$, then show how to transfer the result to the case where $Y(0)=0$. We start with:

\subsection{The case $Y(0)=1/2$}
\noindent Having in mind (\ref{upsilon}), (\ref{UU}), (\ref{itoMc}) and (\ref{beta}), one can write 
\begin{eqnarray*}
\Upsilon(r)&=& \frac{1}{4} \int_0^r Y^2(s) (1-Y(s))^{2\kappa-1} \d U(s),\\
&=& \frac{1}{4} \int_0^{U(r)} \left(S_Y^{-1}(\beta(s))\right)^2 
\left(1-S_Y^{-1}(\beta(s))\right)^{2\kappa-1} \d s,\\
&=& \frac{1}{4} \int_{S_Y(0)}^{\infty} \left(S_Y^{-1}(s)\right)^2
\left(1-S_Y^{-1}(s)\right)^{2\kappa-1} L_{U(r)}^s(\beta)\d s,\\
&=& \frac{1}{4} \int_0^1 s (1-s)^{\kappa-2} L_{U(r)}^{S_Y(s)}(\beta)\d s,\\
&=& \frac{r}{4} \int_0^1 s (1-s)^{\kappa-2} L_{U(r)/r^2}^{S_Y(s)/r}(\beta_r)\d s.
\end{eqnarray*}
 
\noindent We have successively used the occupation density formula and a
scaling argument in deriving the last two identities. We begin with

\medskip

%Again, as in the proof of Proposition
%ref{star}, we drop $\beta_r$ in $L_{\tau_{\beta_r (.)}^{.}} (\beta_r)$ througho%t. 

\noindent {\bf The lower bound.} Let $u, \delta>0$ be given. The last identity
above coupled with Proposition \ref{star}, (\ref{truc}) and a scaling leads to
\begin{eqnarray}
\p_{1/2}(\Upsilon(r)>ur)
&\ge& \p\left(\int_{1-\nu}^1 s(1-s)^{\kappa-2}
  L_{\tau(\lambda_{-})}^{S_Y(s)/r} ds>4u\right)-\p_{1/2}(\Xi_{r}^c),\nonumber
\\
&\ge& \p\left(\int_{\frac{2}{\lambda_{-} r\nu^{\kappa}}}^{\infty}
 \frac{L_{\tau(1)}^y}{y^{2-1/\kappa}} dy > u_1 r^{1-1/\kappa}\right)
-\p_{1/2}(\Xi_{r}^c),\label{oups}
\end{eqnarray}

\noindent where $0<\nu<1$ and 
$$u_1=u_1( \delta, \kappa, \nu, u)=\frac{2^{4-1/\kappa}  \kappa^{2-1/\kappa} {\lambda_{-}}^{1/\kappa}} {(1-\nu)^2} \: u.$$

\smallskip

\noindent Now, let $l$ be a real number such that $1<\kappa<l$. The Cauchy-Schwarz inequality gives
$$
\left(\int_{\frac{2}{\lambda_{-} r\nu^{\kappa}}}^{\infty}
  \frac{L_{\tau(1)}^y}{y^{2-l/\kappa}} \d y\right)^2  \le 
\left(\int_{\frac{2}{\lambda_{-}r\nu^{\kappa}}}^{\infty} \frac{L_{\tau(1)}^y}{y^{2-1/\kappa}} \d y\right)
\left(\int_{0}^{\infty} \frac{L_{\tau(1)}^y}{y^{2-(2l-1)/\kappa}} \d y\right).$$

\noindent Thus, setting $q\;\ddef\; 1-1/{\kappa}$, the first
probability in (\ref{oups}) is, for all $\eta>0$, 
\begin{equation}
\ge\p\left(\int_{\frac{2}{r \lambda_{-}\nu^{\kappa}}}^{\infty} \frac{
L_{\tau(1)}^y}{y^{2-l/\kappa}} \d y > \sqrt{u_1}
r^{lq+\frac{2l-1}{2}\eta}\right)-\p\left( \int_{0}^{\infty} \frac{L_{\tau(1)}^y}{y^{2-(2l-1)/\kappa}} \d y >r^{(2l-1)(q+\eta)}\right).\label{helo}\\
%&\ge& \p\left({\cal S}_{\kappa/n}>(u^{'}+\eta) r^{nq+\frac{2n-1}{2}\eta}\right)%-\p\left(\int_0^{ \frac{1}{r \nu^{\kappa}} }\frac{L_{\tau(1)}^y}{y^{2-n/\kappa}} %\d y > \eta  r^{nq+\frac{2n-1}{2}\eta}\right).\label{helo}
\end{equation}

\noindent By virtue of (\ref{bessstab1}), the second probability in 
(\ref{helo}) involves a stable random variable of parameter $0<\kappa/(2l-1)<1$ and
hence is equivalent to $r^{1-\kappa-\eta \kappa}$. Indeed, from \cite{Bingham}, p 347, we have that for ${\cal S}_{\alpha}$ a
stable random variable of index $0<\alpha<1$, then for $x$ large enough,
$\p({\cal S}_{\alpha} >x)$ is of order $x^{-\alpha}$. (We say that $u(x)$
is {\it of order} $v(x)$ as $x$ tends to
infinity when $\lim_{x \rightarrow \infty}
u(x)/v(x)$ equals some finite nonzero constant.) 

\medskip

\noindent On the other hand, the first one is, for all
$\epsilon>0$, 
\begin{eqnarray}
&\ge& \p\left(\int_{0}^{\infty} \frac{
L_{\tau(1)}^y}{y^{2-l/\kappa}} \d y > (\sqrt{u_1}+\epsilon)
r^{lq+\frac{2l-1}{2}\eta}\right)
-\p\left(\int_{0}^{\frac{2}{\lambda_{-}r\nu^{\kappa}}} \frac{
L_{\tau(1)}^y}{y^{2-l/\kappa}} \d y > \epsilon
r^{lq+\frac{2l-1}{2}\eta}\right),\nonumber \\
&\ge& \p\left(\int_{0}^{\infty} \frac{
L_{\tau(1)}^y}{y^{2-l/\kappa}} \d y > (\sqrt{u_1}+\epsilon) r^{lq+\frac{2l-1}{2}\eta}\right)-\p\left(\sup_{y\ge0}L_{\tau(1)}^{y}>\epsilon \; \nu^{l-\kappa}\; r^{l-1+\frac{2l-1}{2}\eta}\right).\label{hello}
\end{eqnarray}

\noindent Once again, since the first
probability in (\ref{hello}) involves a stable variable of parameter
$\kappa/l<1$, it is of order $r^{1-\kappa-\kappa \eta (1-1/2l)}$.

\smallskip

 Finally, choosing $\nu=r^{-1+(\eta \kappa)/ (l-\kappa)}$ and making use of (\ref{bess1}) tell us that the second probability in (\ref{hello}) is equal to $\epsilon^{-1}
r^{1-\kappa-\eta(\kappa+(2l-1)/2)}$. Putting all that together and having in
mind (\ref{probastar}), and the fact that $u_1$ approaches a constant $u_1 (\delta,
\kappa, o, u)=2^{4-1/\kappa} \kappa^{2-1/\kappa} {\lambda_{-}}^{1/\kappa}$, as
$r$ goes to infinity, we see that $\p_{1/2}(\Upsilon(r)>ur)$ is bounded from
below by some constant times $r^{1-\kappa-\kappa \eta (1-1/2l)}$, for $r$ big 
enough, with $\eta, \epsilon>0$ as small as desired. This completes the proof of (\ref{upsilon:lb}). We now move to  

\bigskip

\noindent {\bf The upper bound.} 
Let $u>u^0 \stackrel{def}{=}\kappa+1/ (\kappa(\kappa-1))$ be given. There exists $\delta_0=\delta_0 (u)>0$
such that for all $0<\delta<\delta_0$,
$u>u_{+}^0\stackrel{def}{=}u^0 (1+\delta)$.
Hence, for such a $\delta$, by virtue of (\ref{ineg}),
$$
\p_{1/2}\left(\Upsilon(r)>ur\right)\le\p\left( \int_0^1 x(1-x)^{\kappa-2}
  L_{\tau(\lambda_{+})}^{S_Y(x)/r} \: \d x>4u\right)+\p_{1/2}(\Xi_{r}^c)
.$$

\noindent Just as in the proof of Proposition \ref{star},
with the same $\mu$ and $\nu>0$ as before, namely $\mu=e^{-r^{\theta}}$ and $\nu=
r^{-\theta/{\kappa}}$, for $0<\theta<1$, we have 
$$\p_{1/2}\left(\Upsilon(r)>ur\right)\le {\mathrm{(I})}+{\mathrm{(II)}}+
{\mathrm{(III)}}+ \p_{1/2}(\Xi_{r}^c),$$

\noindent where
\begin{eqnarray}
{\mathrm {(I)}}&\stackrel{def}{=}& \p\left(\sup_{x\ge 0} L_{\tau(1)}^x >
  \frac{4(u-u_{+}^0)(1-\mu)^{(2-\kappa)^{+}}}{\lambda_+
  {\mu^2}}\right),\nonumber \\
{\mathrm {(II)}}&\stackrel{def}{=}& 2 \p\left( \sup_{0 \le x \le 
\frac{S_Y (1-\nu) \vee |S_Y(\mu)|}{\lambda_{+}r}}
    L_{\tau(1)}^x  >
  \frac{4 u_{+}^0 + (u-u_{+}^0)}{\lambda_{+}\int_{\mu}^{1-\nu} x (1-x)^{\kappa-2}
  dx}\right),\nonumber \\
{\mathrm {(III)}}&\stackrel{def}{=}&\p\left( \int_{g r^{\theta-1}}^{\infty}
  \frac{L_{\tau(1)}^y}
{y^{2-1/\kappa}} \d y >u_2 r^q\right).\label{III}
\end{eqnarray} 

\noindent where $g =(2\lambda_{+} \kappa)^{-1}$ and $u_2=2^{-2+1/\kappa}
  \lambda_{+}^{-1/\kappa}(u-u_{+}^0)$. 

\medskip

We know from (\ref{bess1}) that ${\mathrm {(I)}}$ is of order 
$\mu^2 = e^{-2 r^{\theta}}$, given the choice of $\mu$. On the other hand, for $r$ big 
enough, $\lambda_{+}\int_{\mu}^{1-\nu} x (1-x)^{\kappa-2}  dx$ approaches $4 u_{+}^0$ in which case $(4 u_{+}^0 + (u-u_{+}^0))/4 u_{+}^0$ is greater than 1, implying that ${\mathrm {(II)}}$ decays exponentially fast to zero, by virtue of (\ref{bess2}).

\medskip

So, keeping in mind (\ref{star}), all we need to prove is that 
\begin{equation}\limsup_{r \rightarrow \infty}\frac{\log {\mathrm {(III)}}}{\log r}\leq 1-\kappa.\label{bom}\end{equation}

\noindent (Note that ${\mathrm {(III)}}$ gave us the right order for the lower bound.) The strategy here is akin to the one we used for bounding $\p_{1/2}(J_3 (r, \nu))>\delta)$ from above in the proof of (\ref{probastar}). 

\medskip

\noindent {\bf Second use of iteration:}

\smallskip

For all $\kappa>1$, there exists an integer $n=n(\kappa) \ge 1$ such 
that $2^{n-1} < \kappa \leq 2^n$. Let us suppose first that $2^{n-1}<\kappa
<2^n$. Thus, iterating $n$ times integration by parts followed by a time change, exactly as in the proof of (\ref{probastar}), leads to
\begin{equation}\label{bonn}
{\mathrm {(III)}}\le n e^{-g_n r^{1-\theta/\kappa}}+n e^{-{\hat g}_n \log^2 r}+ \p\left(
    \int_{gr^{\theta-1}}^{\infty}\frac{L_{\tau(1)}^y}{y^{2-2^n /\kappa}} dy
    >\frac{r^{2^{n}q}}{(\log r)^{2^{n+1}-2}}\right),
\end{equation}
\noindent where $g_n$ and ${\hat g}_n$ are two positive real numbers which do {\it not}
    depend on $r$. Note that $1-\theta/\kappa >0$, since $\theta<1<\kappa$.

\medskip
 
Now, as $\kappa / 2^n < 1$, (\ref{bessstab1}) applies, and thus the last
probability is less than or equal to the probability that a stable variable of
index $\kappa / 2^n$ be greater than $r^{2^n q}\: (\log r)^{2-2^{n+1}}$. And
(\ref{bom}) then follows.

\medskip

For $\kappa= 2^n$, in which case stable variables are of no help, all one
needs to prove, given the previous reasoning, is
$$\limsup_{r \rightarrow \infty} \frac{1}{\log r}\; \log \p\left(   
    \int_{gr^{\theta-1}}^{\infty}\frac{L_{\tau(1)}^y}{y} \: dy
    >\frac{r^{\kappa-1}}{(\log r)^{2\kappa-2}}\right) \leq 1-\kappa.$$

\noindent With the help of (\ref{cauchy}), the last probability is, for large
    $r$,  
\begin{equation}
\leq \p\left( {\cal C}_1 > \frac{r^{\kappa-1}}{(\log r)^{2\kappa-1}}\right)
+ \p\left(\int_0^{gr^{\theta-1}} \frac{L_{\tau(1)}^y -1}{y} \: dy < -\frac{r^{\kappa-1}}{2(\log r)^{2\kappa-2}}\right),\label{alou}\end{equation}

\smallskip

\noindent with ${\cal C}_1$ denoting a completely asymmetric Cauchy variable of
index 1. We know that
$$\limsup_{r\rightarrow \infty} \frac{1}{\log r} \: \log  \p\left( {\cal C}_1 > \frac{r^{\kappa-1}}{(\log r)^{2\kappa-1}}\right)\leq 1-\kappa.$$

Hence, we are to prove the same result for the second probability in
(\ref{alou}). To this end, It\^o's formula for $\log u \times (1-L_{\tau (1)}^u)$ reads
$$
\log (g r^{\theta-1})\times (1-L_{\tau (1)}^{g r^{\theta-1}})=
\int_0^{g r^{\theta-1}} \frac{1-L_{\tau (1)}^u}{u}\: du -\;2\; \int_0^{g r^{\theta-1}} \log u \: \sqrt{L_{\tau (1)}^u} \;  d\gamma(u),
$$ 

\noindent with probability one, where we have used the fact that
$$\lim_{u \rightarrow 0} \frac{1-L_{\tau (1)}^u}{u^{\psi}} \; u^{\psi}
\log u = 0, \: \; \; {\mathrm \p-a.s.}$$

\noindent for any $\psi > 1/2$. As a result, the second probability
involved in (\ref{alou}) is, for $r$ large enough,
$$\leq \p\left(L_{\tau (1)}^{g r^{\theta-1}} > \frac{r^{\kappa-1}}
{(\log r)^{2 \kappa}}\right)+\p\left( \int_0^{g r^{\theta-1}} (\log u)^2 \: L_{\tau
  (1)}^u \;  du > \frac{r^{2\kappa-2}}{(\log r)^{4\kappa-2}}\right)+ 2 \: 
\ee^{-\log^2 r}.
$$
\noindent As is easily verified, the first probability above decays
exponentially fast to zero as $r$ approaches infinity. Now Chebychev's inequality together with the fact that
$\int_0^x \log^2 u \: du = x^2 \log^2 x -2x\log x+2x$ implies that the second probability is
$o(r^{1-\kappa})$. So (\ref{bom}) follows, matching the claim.

\medskip

\subsection{The case $Y(0)=0$}

\noindent Let $T_{1/2}=\inf\{s: Y(s)=1/2\}$. We shall need the following
result:

\begin{lemma}\label{l:quick}
For all $n\ge 0$, we have that 
\[
\e_0 (T_{1/2}^n)<\infty.\]
\end{lemma}

\noindent We postpone the proof of this Lemma \ref{l:quick} to the end of the
subsection, first showing how it will be applied. Keeping the same notation as before, we have:

\medskip
\noindent {\bf The upper bound:} Let $v>(1+\kappa)/(\kappa(\kappa-1))$. For $0<\epsilon<1$, 
we have 
\begin{eqnarray}
\p_0(\Upsilon(r)>vr)&=&\p_0 \left( \int_0^r \frac{Y(u)}{(1-Y(u))^2} \d u > vr\right)\le\p_0 \left( 
\int_0^{T_{1/2} +r} \frac{Y(u)}{(1-Y(u))^2} \d u > vr\right), \nonumber \\
&\le& \p_0 \left(\int_0^{T_{1/2}} \frac{Y(u)}{(1-Y(u))^2} \d u > \epsilon vr\right)+\p_0 \left(\int_{T_{1/2}}^{T_{1/2} +r} \frac{Y(u)}{(1-Y(u))^2} \d u > (1-\epsilon)vr\right),\nonumber\\
&\le&\p_0 (T_{1/2} >\epsilon v r/2)+ \p_{1/2} \left(\int_0^r \frac{Y(u)}{(1-Y(u))^2} \d u > 
(1-\epsilon)vr\right),\label{tata}
\end{eqnarray}
\noindent 
where we have used the fact that $T_{1/2}$ is $\p_0$-almost surely finite
(provided by Lemma \ref{l:quick}) together with the strong Markov property. 

\medskip
\noindent On the other hand, by the same reasoning we have:

\bigskip

\noindent {\bf The lower bound:} For all $v>0$ and all $\epsilon$ such that $0<\epsilon<1$,
\begin{eqnarray}
\p_0 \left(\int_0^r \frac{Y(u)}{(1-Y(u))^2} \d u > vr\right)&\ge& \p_0 \left(\int_{T_{1/2}}^r \frac{Y(u)}{(1-Y(u))^2} \d u > vr; \; T_{1/2} \le (1-\epsilon) r\right),\nonumber \\
&\ge& \p_0 \left(\int_{T_{1/2}}^{T_{1/2} + \epsilon r} \frac{Y(u)}{(1-Y(u))^2} \d u > vr; \; T_{1/2} \le (1-\epsilon) r\right),\nonumber \\
&\ge& \p_{1/2} \left(\int_0^{\epsilon r} \frac{Y(u)}{(1-Y(u))^2} \d u > vr\right)-\p_0 (T_{1/2} > 
(1-\epsilon)r).\label{tonton}\end{eqnarray}

\noindent 

\medskip

Now, Markov's inequality in conjunction with Lemma \ref{l:quick} tells us that for all $w>0$,
\[
\limsup_{r\rightarrow \infty}\frac{1}{\log r} \log \p_0 (T_{1/2}>wr)=-\infty,\]

\noindent since this limsup is less than or equal to $-n$ for any integer $n$; we
are done by sending
$n$ to infinity. 

\medskip
Accordingly, choosing $\epsilon$ so small that $(1-\epsilon)w>
(1+\kappa)/(\kappa(\kappa-1))$, since we have proved that (\ref{upsilon:ub}) and (\ref{upsilon:lb}) hold true when $Y(0)=1/2$,
(\ref{tata}) and (\ref{tonton}) deliver (\ref{upsilon:ub}) and
(\ref{upsilon:lb}), as desired.

\bigskip

\noindent Now we turn to Lemma \ref{l:quick}. For completeness, we actually give two proofs, for there are two different ways of writing the Laplace transform of $T_{1/2}$, starting from 0. In the first proof, we exploit the fact that 0 is an entrance boundary. We begin with the Laplace transform of the first exit time of the interval $[l,1/2]$ starting 
from $l<x<1/2$, find its moments and then first send $l$ then $x$ to zero. In the 
second approach, we express the Laplace transform in terms of the hypergeometric 
function, then use results for special functions. Each has its advantages: while the 
first proof provides the finiteness of the moments by 
induction, the second proof, though technically much heavier, gives an explicit formula for 
the moments of $T_{1/2}$. 

\medskip

\noindent {\bf First proof of Lemma \ref{l:quick}:} Recall that 0 is an entrance boundary, unattainable from ]0,1[. Hence, we may write
\begin{equation}\label{shh}
\e_0 (T^n_{1/2})=\lim_{x\rightarrow 0} \lim_{l\rightarrow 0} \e_x ((T_l \wedge T_{1/2})^n),\end{equation}

\noindent for all $n\ge 0$ and $l<x<1/2$. It is well-known that the Laplace transform of $T_l \wedge T_{1/2}$ (=$\inf(T_l, T_{1/2}))$, 
\[
u^l (x)\stackrel{def}{=}\e_x \left(\exp(-\lambda T_l \wedge T_{1/2})  \right), \;\; \;\lambda>0, \]

\noindent satisfies ${\cal L}_Y u^l (x) = \lambda u^l (x),$ $l<x<1/2$ where ${\cal L}_Y$ is the 
infinitesimal generator of $Y$. See for instance \cite{K-T} pages 196-197. Setting 
\[u_n^l (x)\stackrel{def}{=}\e_x ( (T_l \wedge T_{1/2})^n),\]

\noindent we have $u_0^l  \equiv 1$ and
\[{\cal L}_Y u_n^l (x)= -n u_{n-1}^l (x), \;\; \; l<x<1/2,\]       

\noindent with the boundary conditions $u_n^l (l)=u_n^l (1/2)=0,$ for all $n>0$. Making $g\equiv nu_{n-1}^l$ in display (3.11), p 197, \cite{K-T}, leads to  
\begin{eqnarray*}
\frac{2}{n}  u_n^l (x)= \frac{S_Y (l,x)}{S_Y (l,1/2)} \int_x^{1/2} S_Y (t,1/2) (1-t)^{\kappa} u_{n-1}^l (t) \d t \\
 + S_Y (x,1/2) \int_l^x  \frac{S_Y (l,t)}{S_Y (l,1/2)} (1-t)^{\kappa} u_{n-1}^l (t) \d t, \end{eqnarray*}

\noindent for all $l<x<1/2$, $S_Y (a,b)$ denoting $S_Y (b)-S_Y (a)$, for $0<a,b<1$. Recall from (\ref{SY}) the definition 
of $S_Y$. 

\medskip We would like to show that $u_n (0)\stackrel{def}{=} \e_0 (T_{1/2}^n)<\infty.$ We shall do so by induction. Suppose that $u_{n-1} (0)<\infty$, for $n>1$. (This trivially
holds for $n=1$.)
%Now, making $n=1$, so that $u_0 \equiv 1$, letting $l$ tend to 0, at fixed $x$, then letting $x$ tend to 0, we get that 
%\[\e_0 (T_{1/2}= \lim_{x\rightarrow \infty} \e_x (T_{1/2})= \lim_{l \rightarrow 0} u^l_1 (x)= 1/2 \; \int_0^{1/2} 
%S_Y (t,1/2) (1-t)^{\kappa} \d t<\infty.\]
\medskip

\noindent Noting that $u_{n-1}^l (t)\le u_{n-1} (0)$, for all $l<t<1/2$, the
second term in the sum above is less than or equal to $S_Y(x,1/2)\cdot u_{n-1}
(0)\cdot \int_0^x (1-t)^{\kappa} \d t$; this is finite, independent of $l$ and tends to zero as $x$ goes to 0. 

\medskip

Moreover, at fixed $x$, $S_Y (l,x)/S_Y (l,1/2)$ approaches 1 as $l$ tends to 0 since in this case $S_Y (l)$ is 
of order $\log l$. Therefore, sending first $l$ to 0, at fixed $x$, then $x$ to 0, we have
\[\e_0 (T^n_{1/2})= \lim_{x\rightarrow 0} \e_x (T^n_{1/2})= \lim_{x \rightarrow 0} \lim_{l \rightarrow 0} u^l_n (x)= \frac{n}{2} \; \int_0^{1/2} 
S_Y (t,1/2) (1-t)^{\kappa} u_{n-1} (t) \d t<\infty,\]

\noindent by the monotone convergence theorem. We have proved that $u_n
(0)<\infty$, as desired. 
%in other words that $T_{1/2}$ has moments of all order.  

\bigskip

\noindent {\bf Second proof of Lemma \ref{l:quick}:} Here we write the Laplace transform of $T_{1/2}$ in a different way. Let $F(a,b,c,x)$ be the hypergeometric function. (See e.g. \cite{Special}.) The function $F$ solves the following Gaussian differential equation:
\[x(1-x) y^{''} (x)+(c-(a+b+1)x)y^{'}(x)=ab\; y(x).\]

\noindent In the light of (\ref{jacobieds}) this can be rewritten as 
\[{\cal L}_Y y (x)=2ab\; y(x),\]

\noindent for $c=1$ and $a,b>0$ such that $a+b=1+\kappa$. Therefore $(F(a,b,1,Y(t)) e^{-2abt})_{t\geq 0}$ is a local martingale. 
We apply the optional stopping theorem, getting 
\[\e_0 (e^{-2\theta T_{1/2}})=\frac{F(a,b,1,0)}{F(a,b,1,1/2)}=\frac{1}{G(\theta)},\]

\noindent where $\theta=\theta(a)\stackrel{def}{=}ab>0$ and 
\[
G(\theta)\stackrel{def}{=}F(a,b,1,1/2)=\sum_{n\ge 0} \frac{\Gamma(a+n)\Gamma(b+n)}{\Gamma(a) \Gamma(b)(n!)^2} \;
 \frac{1}{2^n}. \]

\noindent It follows that 
\[\e_0 (T_{1/2}^n)=(-1/2)^n \partial_{\theta}^n  (1/G(\theta)\big|_{\theta=0}.\]

\noindent We note that sending for instance $a$ to 0 and $b$ to $1+\kappa$ sends
 $\theta$ to 0. 

\medskip

We begin by rewriting  
$G(\theta)$ as:
\[G(\theta)=1+\sum_{n\ge 1} P_n (\theta) \frac{1}{2^n (n!)^2},\]

\noindent where 
\[P_n (\theta) = a(a+1)...(a+n-1) b(b+1)...(b+n-1)= \prod_{0\le i \le n-1} (\theta+(1+\kappa)i+i^2).\]

We shall show that $G^{'}(0)=2\; \e_0 (T_{1/2} )<\infty$. For higher order
derivatives, the proof follows the same pattern (though with heavier expressions!). 

\medskip

\noindent To this end, we compute the logarithmic derivative of $P_n (\theta)$: 
\[P^{'}_n (\theta) =P_n (\theta) \; \sum_{0\le i\le n-1}
\frac{1}{\theta+(1+\kappa)i+i^2}.\]

\noindent It follows that for all $\theta<\delta$, with $\delta>0$ fixed,  since  both $P_n (\theta)$ and $P_n (\theta)/\theta$ are
increasing functions, $P^{'}_n (\theta)/2^n (n!)^2$ is  bounded from above by $P_n
(\delta) (1/ \delta+ \pi^2/6)/2^n (n!)^2,$ which is summable. By  Lebesgue's dominated convergence theorem followed by the monotone convergence theorem, we arrive at:
 \begin{eqnarray*}
2 \; \e_0 (T_{1/2} )= G^{'}(0)&=& \sum_{n\geq 1} P^{'}_n (0) \; \frac{1}{2^n (n!)^2},\\
&=& \sum_{n\ge 1} \frac{\Gamma(1+n+\kappa)}{\Gamma(2+\kappa)} \; \frac{1}{2^n n \; n!}<\infty,\end{eqnarray*}

\smallskip

\noindent as announced. The next task is to provide the

\section{End of the proofs of Theorems \ref{t:H} and \ref{t:main}.}
\label{s:End}

\noindent Having proved (\ref{upsilon:ub}) and (\ref{upsilon:lb}), we are done with Theorem
\ref{p:derniere}, and hence with Theorem \ref{t:H}, so long as we prove (\ref{I_1-tail}). This is the aim of the following subsection. 

\subsection{Proof of (\ref{I_1-tail})}
\label{half}
There remains to prove:
$$
\limsup_{r \rightarrow \infty} \frac{1}{\log r} \log \p(I_1 (\infty)>wr)\le 1-\kappa,$$

\noindent for all $w>0$. The strategy is akin to the one we used in the previous section for
the tail estimates of $T_{1/2}$; we need only check:

\begin{lemma}\label{ouf}
For any $\alpha$ with $1<\alpha<\kappa$, 
\begin{equation}
\label{nono} \e(I_1^{\alpha} (\infty))<\infty.\end{equation}
\end{lemma}

\smallskip

\bigskip

\noindent {\bf Proof of Lemma \ref{ouf}:} We go back to the identity in law
provided by (\ref{I_1-law}) with $r=\infty$, which makes sense according to (\ref{I_1fini}), and begin with conditioning upon $\{ W(y)\}_{y\in \r},$ and $\xi$, so that the only randomness in the right-hand side of (\ref{I_1-law}) comes from the 0-dimensional squared Bessel
process $Z$. We write $\e^{W,\xi} (\cdot)\; \ddef \; \e( \; \cdot \; | \,
\{ W(y)\}_{y} , \xi)$ for brevity. H\"older's
inequality, for any $\alpha\in (1, \kappa)$ tells us that 
$$
{\e^{W,\xi}} \left( I_1^{\alpha}(\infty) \right) \le S^{\alpha}(\infty) \xi^{\alpha}\,
\left( \, \int_{-\infty}^0 \ee^{-\beta W(y)} g^{\beta}(y) \d y
\right)^{\alpha/\beta} \, { \e^{W,\xi}}\left(
\int_{-\infty}^0 g^{-\alpha}(y)\, Z^{\alpha}\left( { |S(y)| \over S(\infty) \,
\xi}\right) \d y\right) ,
$$

\noindent where $g(y)\; \ddef \; \ee^{W(y)/\alpha}(|S(y)|^{\alpha-1}+
S^{\alpha-1}(\infty)\xi^{\alpha-1})^{1/(\alpha \beta)}$, and $\beta>1$ is such that $\beta^{-1}+\alpha^{-1}=1$.

\medskip

Since $Z$ is a 0-dimensional squared Bessel process starting from
1, we can estimate its moments via its semi-group, see \cite{R-Y} page 441. Indeed, for any $b>0$, 
$$\e(Z^b(t))=\frac{1}{2t}\, e^{-1/2t}\int_0^{\infty} x^{b-1/2} e^{-x/2t} I_1 (\sqrt{x}/t) \d x,$$
\noindent where $I_1$ is the modified Bessel function of index 1, see e.g. \cite{R-Y} p. 549. Plugging the expression for $I_1$ into the above integral and using a Fubini-Tonelli argument followed by the change of variables $y=x/(2t)$, we have 
$$t^{1-b} \, \e(Z^b(t)) = 2^b e^{-1/2t} \sum_{n\geq 0} \frac{\Gamma(n+b+1)}{n!(n+1)!} 
\left(\frac{2}{t}\right)^n.$$

\noindent It is then easily checked that $\limsup_{t \rightarrow \infty} t^{1-b} \, \e(Z^b(t))<\infty$. This implies that there exists $h_{1}(b)$ such that 
$$
\e(Z^b(t))\le h_{1}(b) \, (1+t^{b-1}), \qquad t\ge 0.
$$

\noindent As a consequence, taking $b=\alpha$, we obtain 
$$
{\e^{W,\xi}} \left( \, \int_{-\infty}^0 g^{-\alpha}(y)\, Z^{\alpha}\left( {
|S(y)| \over S(\infty) \, \xi}\right) \d y\right) \le h_{1}(\alpha) \,
\int_{-\infty}^0 g^{-\alpha}(y)\, \left( 1+ { |S(y)|^{\alpha-1} \over
S^{\alpha-1}(\infty) \, \xi^{\alpha-1}} \right) \d y. 
$$

\noindent
$$
{\e^{W,\xi}} \left( I_1^{\alpha}(\infty) \right) \le h_{1}(\alpha) \, S(\infty)
\xi \, \left( \, \int_{-\infty}^0 \ee^{-W(y)} (|S(y)|^{\alpha-1}+
S^{\alpha-1}(\infty)\xi^{\alpha-1})^{1/\alpha} \d y\right)^{\alpha}.
$$

\noindent Making use of  (\ref{triangle}), this leads to
\begin{eqnarray*}
    {\e^{W,\xi}} \left( I_1^{\alpha}(\infty) \right) 
 &\le& h_{2}(\alpha) \, S(\infty) \xi \, \left( \, \int_{-\infty}^0
    \ee^{-W(y)} (|S(y)|^{1-1/\alpha}+ S^{1-1/\alpha}(\infty)\xi^{1-1/\alpha}) \d
    y\right)^{\alpha}
\cr &\le& h_{3}(\alpha) \, S(\infty) \xi \left( \, \int_{-\infty}^0
    \ee^{-W(y)} |S(y)|^{1-1/\alpha} \d y\right)^{\alpha}  
\cr && \qquad \qquad + h_{3}(\alpha) \,S^{\alpha}(\infty) \xi^{\alpha} \,
    \left( \, \int_{-\infty}^0 \ee^{-W(y)} \d y\right)^{\alpha} ,
\end{eqnarray*}

\noindent with $h_{2}(\alpha) =h_{1}(\alpha) d_1^{1/\beta}(\alpha-1)$ and
$h_{3}(\alpha) =h_{2}(\alpha) d_1(\alpha)$. We now take the expectation on both
sides. First, since $\xi$ is exponential of mean $2$, it has
finite moments of all orders. On the other hand, thanks to (\ref{Gamma}), $\e[S^b(\infty)]<\infty$
whenever $b<\kappa$. Moreover, for the same reason as before, since $\alpha<\kappa$, we have 
that $\e [ \, \int_{-\infty}^0 \ee^{-W(y)} \d y ]^{\alpha}<\infty$. Accordingly, for any $\alpha\in (1, \kappa)$, we have:
$$
\e \left( I_1^{\alpha}(\infty) \right) \le h_{4}(\alpha,\kappa) \, \e\left(
\, \int_{-\infty}^0 \ee^{-W(y)} |S(y)|^{1-1/\alpha} \d y\right)^{\alpha}
+ h_{5}(\alpha,\kappa).
$$

It remains for us to handle the expectation term on the right-hand side.
By Lamperti's representation (\ref{Lamperti}), we have
$\int_{-\infty}^0 \ee^{-W(y)} |S(y)|^{1-1/\alpha} \d y \; \law \;
16\int_0^\infty t^{1-1/\alpha} R^{-4}(t) \d t$, where $R$ is a Bessel
process of dimension $(2+2\kappa)$ starting from $R(0)=2$. 

\smallskip
By Lemma \ref{l:Bessel-ub}, this yields $h_{6}(\alpha)\; \ddef \;\e(\,
\int_{-\infty}^0 \ee^{-W(y)} |S(y)|^{1-1/\alpha} \d y)^{\alpha}<\infty$. As a
consequence, for any $\alpha \in (1, \kappa)$, 
$$
\e \left( I_1^{\alpha}(\infty) \right) \le h_{4}(\alpha,\kappa) \,
h_{6}(\alpha) + h_{5}(\alpha ,\kappa),
$$

\noindent finishing the proof of Lemma \ref{ouf}. So
(\ref{I_1-tail}) is proved.

\bigskip

We are done with the proof of Theorem \ref{t:H}. It remains to see how the 
'{\it natural duality}' between $H$ and $X$ enables us to translate Theorem
\ref{t:H} into Theorem \ref{t:main}. The strategy of the proof is akin to that
adopted in \cite{DPZ} for the RWRE case. 

\bigskip

\subsection{End of proof of Theorem \ref{t:main}}
\label{s:main}
We begin with:

\smallskip

\noindent {\bf The upper bound.} Clearly, it suffices to show that for any $v\in (0, v_\kappa)$, 
\begin{equation}
    \limsup_{t\to \infty} {\log \p(X(t)<vt)\over \log t} \le
   1-\kappa.
    \label{ub1}
\end{equation}

\noindent Let $\epsilon>0$ be given. If $X(t)<vt$, then $X$ either stays below
the level $(v+\epsilon)t$ during $[0,t]$, or hits $(v+\epsilon)t$ at time 
$H((v+\epsilon)t)\le t$ and then comes below $vt$ before time $t$. Accordingly,
\begin{eqnarray}\label{bbb} \lefteqn{\p(X(t)<vt)\nonumber}\\
&\leq & \label{freez}
 \p(H((v+\epsilon)t)>t)+\p(H((v+\epsilon)t) \le t ; 
\inf_{s \geq H((v+\epsilon)t)} X(s)<vt).\end{eqnarray}

\noindent Having in mind the definition of the annealed probability $\p$, the second term on the right-hand side is less than or equal to
\begin{equation}E_Q\left(P^W\left(\inf_{s \geq H((v+\epsilon)t)} X(s)-X(H((v+\epsilon)t))
<-\epsilon t\right)\right).\label{ruth}\end{equation}

\noindent Let $\Theta$ be the shift operator, defined by 
$$
\Theta_x W(y)=W(x+y)-W(x).$$

\noindent By virtue of the strong Markov property and the invariance of $\p$ 
under the action of the group $\{\Theta_x,\: x\in \r\}$, the quantity
(\ref{ruth}) equals
\begin{eqnarray*}\lefteqn{E_Q\left(P^W_{(v+\epsilon)t}\left(
\inf_{s\geq0}X(s)-(v+\epsilon)t<-\epsilon t\right)\right)}\\
& & = E_Q\left(P^{{\Theta}_{(v+\epsilon)t}W}\left(
\inf_{s\geq0}X(s)<-\epsilon t\right)\right)\\
& & = \p\left(\sup_{s\geq0}\: (-X(s))\: > \epsilon t\right).\end{eqnarray*}

\noindent Now, thanks to \cite{Kawazu-Tanaka1}, the last probability
approaches zero exponentially fast as $t$ goes to infinity. Accordingly,
taking logarithm of (\ref{bbb}), dividing by $\log r$ then taking the
$\limsup$, and using (\ref{t:H-ub}), since $\epsilon$ is arbitrary, we have the upper bound of $1-\kappa$.\\

\noindent {\bf The lower bound.} Since $G$ is open and separated from
$v_{\kappa}$, it suffices to establish the lower bound for
$G=(v-2\epsilon,v)$, where 
$0<2\epsilon< v< v_{\kappa}$. We set
$${\cal L}_y=\sup_{t \geq H(y)}\left(y-X(t)\right),$$ and 
observe that the event $\{X(t)/t \in (v-2\epsilon,v)\}$ contains the event 
\begin{eqnarray*} & &
  \left\{\frac{(v-2\epsilon)}{v_{\kappa}}t<H((v-\epsilon)t)
<t \; ; 
H(vt)>t ; \; {\cal L}_{(v-\epsilon)t}<\epsilon t\right\}\\
&   &\stackrel{def}{=} A_t \cap B_t \cap C_t.\end{eqnarray*}

\medskip
\noindent Clearly,\begin{eqnarray} \p\left(\frac{X(t)}{t}\in (v-2\epsilon,v)
\right)&\geq & \p(A_t \cap B_t \cap C_t),\nonumber \\
&\geq& \p(B_t|A_t) \p(A_t)- \p(C_t^c). \label{topcool}\end{eqnarray} 

\medskip

\noindent Now, since $\kappa>1$, we know from \cite{Kawazu-Tanaka2} that $H(r)/r$ approaches
$4/(\kappa-1)=v_{\kappa}^{-1},$ $\p$-almost surely, as $r$ tends to
infinity. Thus, as $v<v_{\kappa}$,
 \begin{equation}\label{sai0}\lim_{t\rightarrow\infty}\p(A_t)=1.\end{equation}

\medskip

\noindent On the other hand, once again the strong Markov property together with 
the invariance of $\p$ under $\{\Theta_x,\: x\in \r\}$ and \cite{Kawazu-Tanaka1} imply that 

\begin{equation}\label{sai1}\p(C_t^c)=\p({\cal L}_{(v-\epsilon)t}>\epsilon t)=\p\left(\inf_{s\geq 0}X(s)<-\epsilon t\right)\end{equation} 
 is exponentially small as $t\rightarrow \infty$.

\medskip

\noindent Lastly, since $H((v-\epsilon)t)$ does not depend on 
$\{W(x); x\geq (v- \epsilon)t\}$, it follows by stationarity that
\begin{eqnarray}\p(B_t|A_t)&\geq&
\p\left(H(vt)-H((v-\epsilon) t)>(1-\frac{v-2\epsilon}
{v_{\kappa}})t\: \mid A_t\right),\nonumber\\
%| \frac{v-2\epsilon}{v_{\kappa}}t<H((v-\epsilon)t)<t \right),\nonumber\\
&=&\p\left(H(\epsilon t)>\left(1-\frac{v-2\epsilon}
{v_{\kappa}}\right)t\right).\label{cool}
\end{eqnarray}

\noindent Putting (\ref{topcool}), (\ref{sai0}), (\ref{sai1}), (\ref{cool}) and (\ref{t:H-lb}) together completes the proof of the lower bound in Theorem \ref{t:main}.
\qed

\bigskip

Although we had come up with the iteration scheme as a way of avoiding the
technical difficulty associated to a Sturm-Liouville approach, upon the
prodding of the referee we were in fact able to push through that method as well. So for completeness, we include this approach in the next section. 

\section{A Sturm-Liouville alternative to the iteration scheme}
\label{s:Appendix}

In (\ref{iteration}) and (\ref{III}), the iteration scheme enabled us to prove that 
\begin{equation}
\limsup_{r\rightarrow \infty} \frac{1}{\log r} \; \log \p\left(\int_{\epsilon}^{\infty} \frac{L_{\tau_1}^y}{y^{1+\gamma}} \d y > r^{\gamma}\right)=-\infty,\label{lio}\end{equation}

\noindent 
for $\epsilon=\epsilon(r)=r^{\theta-1}\rightarrow 0$, ($0<\theta<1$), with $\gamma>0$ equal to $1+1/\kappa$
in (\ref{iteration}) and to $q=1-1/\kappa$ in (\ref{III}). 

\bigskip

An alternative way of estimating the tails of $\int_{\epsilon}^{\infty} L_{\tau_1}^y/y^{1+\gamma} \d
y$ is to study its Laplace transform. Thanks to a result of Pitman and Yor \cite{Pitman-Yor}, this reduces to solving 
a Sturm-Liouville equation, as we will see in this section. 

\bigskip

From \cite{Pitman-Yor}, we get that, for all $\lambda>0$, 
\begin{equation}\label{py}
\e\left(\exp\left(-\lambda \int_0^{\infty} \frac{L_{\tau_1}^y}{y^{1+\gamma}} \;
    1_{y\geq \epsilon} \; \d y\right)\right)= e^{{\phi^{'}_{\lambda}} (0^{+})/2},\end{equation}

\noindent with $\phi_{\lambda}^{'} (0^+)/2$ denoting the right-derivative of $\phi_{\lambda}$ at 0, where $\phi_{\lambda}$ is the unique convex, decreasing, nonnegative solution of the following 
Sturm-Liouville equation with $\phi_{\lambda} (0)=1$:
\[\Phi_{\lambda}^{''} (x)-\frac{2\lambda}{x^{1+\gamma}}\; 1_{x\ge\epsilon} \; \Phi_{\lambda} (x)=0.\]

\noindent 
Note that one should a priori multiply $e^{\phi_{\lambda}^{'} (0^{+})/2}$ by
    $\phi_{\lambda} (\infty)^0$ in (\ref{py}). The convention $0^0=1$ allows us to omit this factor.

%Whether or not $\phi_{\lambda} (\infty)>0$ (according 
%as $\int^{\infty}  t^{-\gamma} \d t <\infty$, see \cite{Pitman-Yor}, p 428), $\phi_{\lambda} (\infty)^0=1$ 

\smallskip

%Note that $\phi_{\lambda}^{'}$ is a constant on $[0,\epsilon]$, and that it equals
%$\phi_{\lambda}^{'}(0^{+})=\phi_{\lambda}^{'}(\epsilon)$.

\medskip

Solving the Sturm-Liouville equation amounts to solving the following
Riccati's differential equation with $y=\Phi_{\lambda}^{'} / \Phi_{\lambda}$:
\[y^{'}(x)+y^2(x)= 2\lambda \; x^{-1-\gamma},\; \; \; x\ge \epsilon .\]

\smallskip We find from \cite{Watson}, p 88-89, that this is soluble in
finite terms only when $(1-\gamma)/2$ is the inverse of an odd
integer, that is when $\kappa=n+1/2$, for $n\ge 1$.

\medskip

For arbitrary $\kappa>1$, and $\lambda>0$, the general solution of our Sturm-Liouville equation reads: 
$$\Phi_{\lambda} (x)=\sqrt{x} \; {\cal C}_{\kappa} (i \kappa \sqrt{8\lambda} x^{(1-\gamma)/2}), \; 
\; \; x \ge \epsilon,$$

\noindent with ${\cal C}_{\kappa}$ a cylindrincal function of index
$\kappa$; see \cite{Watson}, pages 82-83. 

\medskip

\noindent  Now $\Phi^{'}_{\lambda}$ is constant on the interval $[0,\epsilon]$; a few lines of computation give that for $\epsilon>0$, 
\[\Phi_{\lambda}^{'}(0^{+}) = \Phi_{\lambda}^{'} (\epsilon) = i\; \sqrt{2\lambda} \epsilon^{-\gamma/2} \; 
\left({\cal C}_{\kappa-1} (i\kappa \sqrt{8\lambda}  \epsilon^{1/(2\kappa)})\; 1_{\gamma=q}\; +\; 
 {\cal C}_{\kappa+1} (i\kappa \sqrt{8\lambda}  \epsilon^{-1/(2\kappa)})\; 1_{\gamma=1+1/\kappa}\;  \right),\]

\smallskip

\noindent with $1_A$ denoting the indicator function of $A$.  
\medskip

\noindent From the analyticity of ${\cal C}_{\kappa}$ one gets that, as a function of $\lambda>0$,  $\Phi_{\lambda}^{'} (\epsilon)$ is analytic, thus for $\lambda>0$ small enough 
(depending on $\epsilon$ or equivalently on $r$), one could write:

\[
2\; \log \; \e\left(\exp\left(\lambda \int_{\epsilon}^{\infty} \frac{L_{\tau_1}^y}{y^{1+\gamma}} \; \d y\right)\right)
= {\phi^{'}_{-\lambda}} (\epsilon)=\sqrt{2\lambda} \epsilon^{-\gamma/2} \; {\cal C}_{\kappa \pm 1} \left(\kappa \sqrt{8\lambda}  \epsilon^{\mp 1/(2\kappa)}\right),\]

\noindent for {\it the} cylindrical function ${\cal C}_{\kappa}$ determined by the particular solution $\phi_{\lambda}$. 

\medskip

A cylindrical function can be expressed as:
\[{\cal C}_{\kappa}(x)=a_{\kappa} \; J_{\kappa} (x) + b_{\kappa} \; Y_{\kappa} (x),\]

\noindent with $a_{\kappa}$ and $b_{\kappa}$ two periodic functions of $\kappa$ with period one, and where $J_{\kappa}$ and $Y_{\kappa}$ are Bessel functions of the first and second kind respectively. From pp 622, 625 and 627 of 
\cite{Morse} (or pp 74 and 199 of \cite{Watson}) we have the asymptotic equivalents 
of $J_{\kappa}$ and $Y_{\kappa}$ at 0 and infinity: for $x$ in the neighborhood of 0, $J_{\kappa}(x)$ is of 
order $x^{\kappa}$ 
and $Y_{\kappa}(x)$ of order $x^{-\kappa}$ (for $\kappa>1$). Furthermore, for $x$ large, both $J_{\kappa}(x)$ and 
$Y_{\kappa}(x)$ are of order $x^{-1/2}$. 

\medskip

This provides all the ingredients for proving our tail estimates. Indeed, for all $u>0$, with the previous choices of 
$0<\theta<1$, $\epsilon=r^{1-\theta}$, and for $\lambda=\lambda(r)$ chosen to go very slowly to zero as r tends to infinity, 
an exponential inequality together with (\ref{py}) yields
\[\p\left(\int_{\epsilon}^{\infty} \frac{L_{\tau_1}^y}{y^{1+\gamma}} \; \d y > u\; r^{\gamma}  \right)
\le \exp(-\lambda u r^{\gamma} + \frac{1}{2}\phi_{-\lambda}^{'} (\epsilon) ). \]

\medskip
\noindent By virtue of the choice of $\lambda$, $\phi_{-\lambda}^{'} (\epsilon)$ is of order 
$\epsilon^{-\gamma}=r^{(1-\theta)q}=o(r^{\gamma})$, for $\gamma=1-1/\kappa$,  
and $r^{(1-\theta)(\gamma+1)/4}=o(r^{\gamma})$, for $\gamma=1+1/\kappa$. 

\medskip

We have proved (\ref{lio}) for both (\ref{iteration}) and (\ref{III}).
\bigskip

\noindent{\bf Acknowledgements} I am indebted to Zhan Shi for a wealth of very
useful
comments, for his unfading support, and for helping me prove inequality (\ref{I_1-tail}). 
Also, warm thanks are due to Yueyun Hu for his encouragement. This work owes a
lot to Alby Fisher for his invitation to the University of S\~ao Paulo, for his
unconditional support and for many fruitful discussions; his passion for
mathematics has been a constant inspiration for me. Finally, thanks are due to the referee for his careful reading, for his constructive suggestions and for pointing out the 
possibility of the Sturm-Liouville alternative to our iteration scheme.

\bigskip

%\end{document}

%% Authors' adresses
%%

{\footnotesize 
\baselineskip=12pt

\noindent 
\begin{tabular}{lll}

    & \hskip30pt Marina Talet \\ 

    & \hskip30pt C.M.I. Universit\'e de Provence \\   

    & \hskip30pt LATP, CNRS-UMR 6632 \\   

    & \hskip30pt 39, rue F. Joliot Curie \\   

    & \hskip30pt F-13453 Marseille Cedex 13 \\  

    & \hskip30pt France \\ 

    & \hskip30pt {\tt marina@cmi.univ-mrs.fr}
\end{tabular}

}%end of authors' adresses

\end{document}